\title{\LARGE \bf%
A Fredholm transformation for the rapid stabilization of a degenerate parabolic equation\footnote{This research was partially supported by the French Grant ANR ODISSE (ANR-19-CE48-0004-01) and
was also conducted in the framework of the regional programme ”Atlanstic 2020, Research, Education
and Innovation in Pays de la Loire", supported by the French Region Pays de la Loire and the European
Regional Development Fund.}
}

\author{Ludovick Gagnon$^{1}$, Pierre Lissy$^{2}$ and Swann Marx$^{3}$}

\documentclass[10pt]{article}

\usepackage{amsfonts,latexsym,amsmath,amssymb}
\usepackage[mathscr]{eucal}
\usepackage{graphicx,color}
\usepackage{comment}
\usepackage{hyperref}
\usepackage{epstopdf}
\newtheorem{theorem}{Theorem}
\newtheorem{lemma}{Lemma}
\newtheorem{proposition}{Proposition}
\newtheorem{corollary}{Corollary}
\newtheorem{definition}{Definition}

\newtheorem{remark}{Remark}

\catcode`\@=11
\def\downparenfill{$\m@th\braceld\leaders\vrule\hfill\bracerd$}
\def\overparen#1{\mathop{\vbox{\ialign{##\crcr\crcr
\noalign{\kern0.4ex}
\downparenfill\crcr\noalign{\kern0.4ex\nointerlineskip}
$\hfil\displaystyle{#1}\hfil$\crcr}}}\limits}
\catcode`\@=12

\def\NN{{\mathbb N}}    
\def\RR{{\mathbb R}}    
\def\CC{{\mathbb C}}    
\def\L{{\cal L}} 

\def\<{\langle}
\def\>{\rangle}

\newcommand{\cqfd}
{%
\mbox{}%
\nolinebreak%
\hfill%
\rule{2mm}{2mm}%
\medbreak%
\par%
}

\begin{document}
\maketitle
\footnotetext[1]{Inria, Villers-l\`es-Nancy, F-54600, France. E-mail : \texttt{ludovick.gagnon@inria.fr}}
\footnotetext[2]{CEREMADE, Universit\'e Paris-Dauphine \& CNRS UMR 7534, Universit\'e PSL, 75016 Paris, France. E-mail: \texttt{lissy@ceremade.dauphine.fr.}}
\footnotetext[3]{LS2N, \'Ecole Centrale de Nantes \& CNRS UMR 6004, Nantes, France. E-mail: \texttt{swann.marx@ls2n.fr}}

\begin{abstract}
This paper deals with the rapid stabilization of a degenerate parabolic equation with a right Dirichlet control. Our strategy consists in applying a \emph{backstepping} strategy, which seeks to find an \emph{invertible} transformation mapping the degenerate parabolic equation to stabilize into an exponentially stable system whose decay rate is known and as large as we desire. The transformation under consideration in this paper is Fredholm. It involves a kernel solving itself another PDE, at least formally. The main goal of the paper is to prove that the Fredholm transformation is well-defined, continuous and invertible in the natural energy space. It allows us to deduce the rapid stabilization. 
\end{abstract}
\\


\section{Introduction}

The topic of this paper is the stabilization of a degenerate parabolic equation with a right Dirichlet control.  Let $\alpha\in(0,1)$. The PDE under investigation is written as follows:

\begin{equation}
\label{system}
\left\{
\begin{aligned}
&\partial_t u=(x^\alpha \partial_x u)_x, \qquad \quad \qquad (t,x)\in (0,T)\times (0,1),  
\\&u(t,0)=0,\: u(t,1)=U(t), \quad t\in (0,T), \\
&u(0,x)=u_0(x), \, \, \quad \quad \quad \quad \quad x\in (0,1),
\end{aligned}
\right.
\end{equation}
where $U$ denotes the control. The initial condition $u_0$ will be chosen in $L^2(0,1)$. A suitable functional setting will be introduced in Section \ref{sec_main} and the well-posedness of \eqref{system} will also be justified  in  Section \ref{ABKS} when the control $U$ is chosen in an appropriate feedback form. Since
$
\alpha\in (0,1),
$
the equation corresponds to a``weakly degenerate'' case in the sense of \cite{CMV05}. \\

In contrast with the null-controllability problem (which has been totally solved for any $\alpha\in(0,1)$ in \cite{CanMarVan,MR2227693} for distributed controls, \textit{i.e.} controls acting in a subdomain included in $(0,1)$, but the well-known extension method also ensures null-controllability for a Dirichlet control at $x=1$), the stabilization problem aims at finding a control depending on the state of the system, leading to a \emph{feedback-law}, so that the origin is exponentially stable. We say that the system is in \emph{closed-loop} when this control is defined in this fashion. One of the main advantage of feedback laws is the robustness of the closed-loop system, \textit{i.e.} the asymptotic stability is preserved when the system is subject to some small disturbances.\\ 

There exist many ways to define such a feedback-law for infinite-dimensional systems. Let us mention for instance \cite{urquiza2005rapid}, where a Gramian approach is followed in order to stabilize rapidly abstract systems with unbounded controls, the computation of Lyapunov functionals for systems of first-order linear hyperbolic PDEs \cite{bastin2016stability} leading then to some feedback-law, or the backstepping technique (see for instance \cite{KrsticSmyshlyaev_Book}), which is the technique that we want to apply to stabilize the origin of \eqref{system}.\\

The backstepping technique consists in finding an integral transformation that maps a (possibly open-loop unstable) system (which will be here \eqref{system}) into an exponentially stable system (which is called the target system). In most cases, this integral transformation gives rise to a kernel solving a PDE. One of the main challenges relies on the analysis of this latter PDE, that is most of the time not a classical one and therefore requires to apply specific techniques. This analysis allows indeed to prove the continuity and the invertibility of the integral transformation, which leads to the exponential stability of the original system, \eqref{system} in our case. \\

The target system that we will consider in this paper is the following :
\begin{equation}
\label{cibleFredholm}
\left\{
\begin{aligned}
&\partial_t v=(x^\alpha \partial_x v)_x-\lambda v, \quad (t,x)\in (0,T)\times (0,1)
\\&v(t,0)=v(t,1)=0, \, \,\, \, \quad t\in (0,T) \\
&v(0,x)=v_0(x), \qquad \quad \, \,\, x\in (0,1).
\end{aligned}
\right.
\end{equation}
The initial condition $v_0$ of this system is exponentially stable. Indeed,  in a suitable functional setting that will be defined later on, the time derivative of $V(t):=\int_0^1 v(t,x)^2 \textrm{dx}$ along the solutions to \eqref{cibleFredholm} yields, for all $t\geq 0$,

$$\begin{aligned}\frac{d}{dt} V(t)&= 2\int_0^1 v(t,x)\partial_t v(t,x)\textrm{dx}\\&= 2\int_0^1 v(t,x)\left((x^\alpha \partial_x v(t,x))_x-\lambda v(t,x) \right)\textrm{dx}\\&= -2\int_0^1 x^\alpha\left(\partial_x v(t,x)\right)^2\textrm{dx}-2\lambda \int_0^1v(t,x)^2\textrm{dx},\end{aligned}$$
whence

\begin{equation}\label{dissv}
\frac{d}{dt} V(t) \leq - 2\lambda V(t).
\end{equation}
Inequality \ref{dissv} implies immediately the exponential stability of the origin with an exponential decay rate of convergence $\lambda$, that we can choose as large as we want. This explains why this kind of stabilization method is referred as a \emph{rapid stabilization}. \\

The most popular integral transformations in the context of the backstepping of PDEs are the Volterra and the Fredholm ones (see \cite{KrsticSmyshlyaev_Book} and \cite{CGM} for the Volterra transformation and the Fredholm one, respectively). To the best of our knowledge, the Volterra transformation has been introduced in \cite{russell1978controllability} and then popularized by Miroslav Krstic and his co-authors (see e.g. \cite{KrsticSmyshlyaev_Book} for a good overview of such a technique). The Fredholm transformation that we will use here has then been applied independently on the Kuramoto-Sivashinky equation, Korteweg-de Vries equation in \cite{CoronLu15} and \cite{CoronLu14}, respectively, and on some hyperbolic PDEs where appear non-local terms in \cite{bribiesca2015backstepping}. Let us notice that this Fredholm transformation allows also to solve the stabilization problem of some PDEs with distributed controls, such as the Schr\"odinger equation \cite{CGM} and the transport equation \cite{MR4001127}. 

We also want to emphasize on the fact that the backstepping technique allows to solve the problem of the null controllability, as illustrated for instance in \cite{coron2017null}, which focuses on a heat equation, or \cite{coron2016stabilization} that is devoted to the case of some hyperbolic PDEs. \\

There exists already a vast literature dealing with the null-controllability of some degenerate parabolic equations similar to \eqref{system}, starting with the seminal work \cite{CMV05}. Dealing with this PDE is not an easy task, because the spatial operator vanishes at $x=0$, which leads to technical difficulties for obtaining results in terms of controllability. Furthermore, this imposes to work in a functional setting that is  different from the usual Sobolev setting for the heat equation (\textit{i.e.} a functional setting involving weighted Sobolev spaces, as we will see later on). The weakly degenerate case that is of interest for us here has been totally solved in \cite{CanMarVan,MR2227693}, where some global parabolic Carleman estimates (similar to the ones proved in \cite{MR1406566} for the usual heat equation) are derived to prove the null-controllability of a parabolic degenerate equation with a distributed control (and as a consequence for Dirichlet boundary controls at $x=1$ by the extension method). The proof of a null-controllability result for a Dirichlet control at  $x=1$ has been solved in \cite{Gue} thanks to the transmutation method.\\

More recently, these results have been improved in \cite{buffe2018spectral}, where a Lebeau-Robbiano strategy (see \cite{MR1312710,MR1620510}) is followed in order to prove the finite-time stabilization of such a system, with distributed  controls. Let us mention that, in contrast with the null-controllability property, our result of rapid stabilization of \eqref{system} with boundary control cannot be deduced from \cite{buffe2018spectral}, since the extension method will not give a feedback law in this context. In the case where the degeneracy occurs at both end-points, a result of distributed controllability has been obtained in \cite{MR2227700}. These results have been generalized in \cite{MR3456387}, based on the flatness approach, with very general coefficients that can notably be weakly degenerate or singular at many points of the interval $(0,1)$. The flatness approach has also been used successfully in \cite{MR3583882} in the case of strong degeneracy (\textit{i.e.} $\alpha\in[1,2]$) and boundary control at $x=1$. For some generalizations in space dimension larger than $1$, we refer to \cite{MR3430764}. \\

The aim of our article is to design a feedback law such that the resulting closed-loop system is exponentially stabilized with an arbitrary decay rate. To the best of our knowledge, this result is new in the context of boundary control of such an equation. This article is organized as follows: Section \ref{sec_main} is devoted to the introduction of some concepts related to parabolic degenerate operator and to the statement of our main results.
\textcolor{black}{Section \ref{sec_wp} gives a a Riesz basis property for some family of functions $\{\psi_n\}_{n\in\mathbb N^*}$ closely related to the kernel of our Fredholm transformation, and some consequences}.  Section \ref{sec_prop} deals with the proof of the continuity and the invertibility of the Fredholm operator that we consider to transform our system \eqref{system} into an exponentially stable one.  Finally, in Section \ref{ABKS}, we prove our main result of rapid exponential stabilization.

\section{Preliminaries and main results}

\label{sec_main}

In this section, we first recollect some results from \cite{CanMarVan, Gue} dealing with the spectrum of degenerate parabolic operators. Second, we introduce the transformation mapping from \eqref{system} to \eqref{cibleFredholm} and the involved kernel to be solved. Third, we state the main results of our paper.


\subsection{Some properties on the degenerate operator $A$}

We recall some well-known facts that can for instance be found in \cite{CanMarVan}.
For $\alpha \in (0,1)$, define 
\[
H^1_\alpha(0,1):=\{ f\in L^2(0,1) \, | \, x^{\alpha/2} f' \in L^2(0,1)\},
\]
endowed with the inner product 
\[
(f,g)_{H^1_\alpha}:=\int_0^1 fg+ x^\alpha f'g' \, \textrm{\textrm{dx}}, \qquad f,g\in H^1_\alpha(0,1).
\]
We call $||\cdot||_{H^1_\alpha}$ the corresponding norm. $(H^1_\alpha(0,1),(\cdot,\cdot)_{H^1_\alpha})$ is a Hilbert space.
We can now define 
\[
H^1_{\alpha,0}(0,1):=\{ f\in H^1_{\alpha}(0,1) \, | \, f(0)=f(1)=0\}.
\]
We recall the following Hardy-Poincar\'e inequality:
\begin{equation}\label{HPI}
\int_0^1 x^\alpha \left| f' \right|^2 \, \textrm{\textrm{dx}} \geq \dfrac{(1-\alpha)^2}{4} \int_0^1 x^{\alpha-2} \left| f \right|^2 \, \textrm{\textrm{dx}}, \qquad \forall f\in C^{\infty}_0(0,1].
\end{equation}
Using the fact that $x^{\alpha-2}\geqslant 1$ for $x\in (0,1)$, we obtain (since $H^1_{\alpha,0}(0,1)$ is the closure of $C^{\infty}_0(0,1)$ endowed with the $H^1_\alpha(0,1)$ norm) 
\begin{equation}\label{HPIi}
\int_0^1 \left| f \right|^2 \, \textrm{\textrm{dx}} \leq \dfrac{4}{(1-\alpha)^2} \int_0^1 x^\alpha |f'|^2 \, \textrm{\textrm{dx}}, \qquad \forall f \in H^1_{\alpha,0}(0,1).
\end{equation}
Finally, we define the following norm on $H^1_{\alpha,0}(0,1)$, which is equivalent to the restriction of the norm $||\cdot||_{H^1_\alpha}$ on $H^1_{\alpha,0}(0,1)$ thanks to \eqref{HPIi}:
\[
\| f \|^2_{H^1_{\alpha,0}}:=\left( \int_0^1 x^\alpha | f' |^2 \, \textrm{\textrm{dx}} \right)^2.
\]
We define the unbounded operator $A : D(A) \subset L^2(0,1) \rightarrow L^2(0,1)$ by:
\[
\begin{cases}
Au:= (x^\alpha u_x)_x, \\ 
D(A):=\{ u\in H^1_{\alpha,0}(0,1) \, | \, x^\alpha u_x \in H^1(0,1) \}.
\end{cases}
\]
The operator $-A$ is self-adjoint, positive definite, and with compact resolvent. Therefore, there exists a Hilbert basis $\{\phi_n\}_{n \in \NN^*}$ of $L^2(0,1)$ and an increasing sequence $(\lambda_n)_{n \in \NN^*}$ of real numbers (than can be proven to be distinct) such that $\lambda_n>0$, $\lambda_n \rightarrow +\infty$ and 
\begin{equation}
\label{ef}
-A\phi_n=\lambda_n \phi_n.
\end{equation}
Before introducing the basis of this operator, let us recall some facts on Bessel functions (for more details on the Bessel functions, we refer to \cite{MR0010746}). For a complex number $\nu$, Bessel functions of order $\nu$ are solutions of the following second order ordinary differential equation:
\begin{equation}\label{edob}
x^2 y^{\prime\prime}(x) + x y^\prime(x) + (x^2-\nu^2)y(x)=0,\: x\in (0,\infty).
\end{equation}
In our analysis, we will focus on the Bessel functions of the first kind, that are finite at the origin (in contrast with the Bessel functions of the second kind that diverge at $x=0$, and are also solution of \eqref{edob}). On can define Bessel functions of the first kind by the following formula:
\begin{equation*}
J_\nu(x)=\sum_{m\geq 0} \frac{(-1)^m}{m!\Gamma(m+\nu+1)}\left(\frac{x}{2}\right)^{2m+\nu},\: x\in\mathbb C,
\end{equation*}
where $\Gamma(\cdot)$ is the gamma function\footnote{Let us recall that the gamma function is defined as follows: $\Gamma:\: z\mapsto \int_0^{+\infty} t^{z-1}e^{-t}dt$}. It can be proved that the function $J_\nu$ has an infinite number of real zeros. These zeros are all simple, with the possible exception of $x=0$, depending on the value of $\nu$.
From now on, we concentrate on the following choice of parameter $\nu$:
\begin{equation}\label{defnu}
\nu := \dfrac{1-\alpha}{2-\alpha}.
\end{equation}
We also introduce the parameter $\kappa$ given by 
\begin{equation}\label{defkappa}
 \kappa:= \dfrac{2-\alpha}{2}.
\end{equation}
We have the following formulas: for all $(a,b)\in\mathbb{R}^2$ such that $a\neq b$, one has
\begin{equation}
\begin{split}
\label{IPB}
\int_0^1 x J_\nu(ax)J_\nu(bx) \textrm{dx} = 
 &\frac{1}{b^2-a^2}\left[b J_{\nu}(a)J^\prime_{\nu}(b) - aJ^\prime_{\nu}(a)J_{\nu}(b)\right],
\end{split}
\end{equation}
and
\begin{equation}
\begin{split}
\label{IPB2}
\int_0^1 x J_\nu(ax)^2 \textrm{dx} = &\frac{1}{2}\left[\left(1-\frac{\nu^2}{a^2} \right) J_{\nu}(a)^2+J_{\nu}'(a)^2\right].
\end{split}
\end{equation}
These two formulas come respectively from \cite[Section 11.2, (2) and (4)]{Yudell}.
We also have the following orthogonality property:
\begin{equation}
\label{PropOrthogonal}
\int_0^1 x^{1-\alpha} J_v(j_{\nu,n}x^\kappa)J_v(j_{\nu,m}x^\kappa) \textrm{dx} = \frac{\delta_{nm}}{2\kappa} \left[J^\prime_\nu(j_{\nu,n})\right]^2,
\end{equation}
where $(j_{\nu,n})_{n\in\mathbb N^*}$ is the sequence of the (positive, simple and ordered in increasing order) zeros of $J_\nu$ and $\delta_{nm}$ is the Kronecker symbol, see \cite[Section 4.3.1]{Gue}.

We are now in position to define the spectrum of the operator $-A$. The normalized eigenvectors $\phi_n$ and eigenvalues $\lambda_n$ introduced in \eqref{ef} are given by 
\begin{equation}
\label{defphi}
\phi_n(x)=\dfrac{(2\kappa)^{1/2}}{J_\nu'(j_{\nu,n})}x^{(1-\alpha)/2}J_{\nu}(j_{\nu,n}x^\kappa), \quad x\in (0,1), \quad n\in \NN^*,
\end{equation}
and 
\begin{equation}\label{dln}
\lambda_n:=(\kappa j_{\nu,n})^2,  \quad n\in \NN^*.
\end{equation}
Let us recall the expansion of the zeros of the Bessel functions for $n\rightarrow\infty$:
\begin{equation}\label{zeroasympt}
j_{\nu,n}=\pi(n+\nu/2-1/4)-\dfrac{4\nu^2-1}{8a}+\mathcal{O}\left(\dfrac{1}{n^3}\right).
\end{equation}
Also recall 
\begin{equation}\label{dern}
 C_1 \sqrt{\dfrac{1}{j_{\nu,n}}}\leq \left|J'_{\nu}(j_{\nu,n})\right| \leq C_2 \sqrt{\dfrac{1}{j_{\nu,n}}}.
\end{equation}

Let us define now the derivative of $\phi_n$ evaluated at $x=1$. It will be indeed useful in our analysis in the sequel of the paper, since it corresponds to the location where is the control. Using the definition of $\phi_n(x)$ given in \eqref{defphi}, we obtain that its derivative with respect to $x$ is as follows:
\begin{equation}
\phi_n^\prime(x) = \frac{(2\kappa)^{1/2}}{J_{\nu}^\prime(j_{\nu,n})}\left[\frac{1-\alpha}{2}x^{-(1+\alpha)/2} J_{\nu}(j_{\nu,n}x^{\kappa}) + x^{(1-\alpha)/2} j_{\nu,n} J^\prime_{\nu}(j_{\nu,n} x^{\kappa})\kappa x^{\frac{-\alpha}{2}}\right].
\end{equation}
Then, evaluating the latter expression at $x=1$ and using the fact that $j_{\nu,n}$ are zeros of the Bessel function $J_\nu$, one obtains
\begin{equation}
\label{behaviorB}
\phi_n^\prime(1) =  (2\kappa)^{1/2} \kappa j_{\nu,n}.
\end{equation}
In all what follows, we will assume the following non-resonance property:
\begin{equation}\label{nonresonance}
\lambda-\lambda_n\not =0,\,\,\lambda-\lambda_n\neq \lambda_k \textrm{ for any } k,n\in \NN^*.
\end{equation}
Notice that is it not a very restrictive condition since the set $\{\lambda_k-\lambda_n|\, k,n\in \NN^*\}$ is discrete, so that $\lambda$ can always be perturbed by some $\varepsilon>0$ as small as we want so that $\lambda+\varepsilon$ verifies \eqref{nonresonance}.

Let us finally remark that we can characterize $D(A)$ and $H^1_{\alpha,0}(0,1)$ as follows: 
$$D(A)=\left\{f=\sum_{n=1}^\infty a_n\phi_n\in L^2(0,1)\left|\sum_{n=1}^\infty \lambda_n^2a_n^2<\infty\right.\right\},$$
so that 
\begin{equation}\label{dap}D(A)'=\left\{f=\sum_{n=1}^\infty a_n\phi_n\in \mathcal D'(0,1)\left|\sum_{n=1}^\infty \frac{a_n^2}{\lambda_n^2}<\infty\right.\right\},\end{equation}
and
$$H^1_{\alpha,0}(0,1)=D(A^{\frac{1}{2}})=\left\{f=\sum_{n=1}^\infty a_n\phi_n\in L^2(0,1)\left|\sum_{n=1}^\infty \lambda_na_n^2<\infty\right.\right\},$$
so that 
\begin{equation}\label{dual12}H^1_{\alpha,0}(0,1)'=D(A^{\frac{1}{2}})'=\left\{f=\sum_{n=1}^\infty a_n\phi_n\in \mathcal D'(0,1)\left|\sum_{n=1}^\infty \frac{a_n^2}{\lambda_n}<\infty\right.\right\}. \end{equation}
All these spaces are endowed with the natural scalar product induced by their definitions. We also emphasize that the operator $A$ can be uniquely extended from $D(A^{1/2})$ to $D(A^{1/2})'$ and from $L^2(0,1)$ to $D(A)'$ (see \cite[Section 3.4]{TucsnakWeissBook}). We still denote by $A$ these extensions, which will be made clear by the context. 

With all these notations, we can already state the following well-posedness result for \eqref{cibleFredholm} by applying general results on semigroup theory  on \textcolor{black}{ $A-\lambda\,I: D(A)\rightarrow L^2(0,1)$ and $A-\lambda\,I: L^2(0,1)\rightarrow D(A')$} (see \textit{e.g.} \cite[Proposition 2.3.5]{TucsnakWeissBook}).
\begin{proposition}\label{rv}
For any  $v_0\in D(A)$, there exists a unique solution $v$ to \eqref{cibleFredholm} verifying moreover 
$$v\in C^1([0,\infty),L^2(0,1)\times C^0([0,\infty), D(A)).$$

For any  $v_0\in L^2(0,1)$, there exists a unique solution $v$ to \eqref{cibleFredholm} (the first line being verified in $D(A)'$ for any $t\in [0,\infty)$) verifying moreover 
$$v\in C^1([0,\infty),D(A)')\times C^0([0,\infty), L^2(0,1)).$$
\end{proposition}

Using Proposition \ref{rv}  and an easy density argument, it is not difficult to deduce the following exponential stability estimate.
\begin{corollary}\label{cdv}
Every solution of \eqref{cibleFredholm} with initial condition $v_0\in L^2(0,1)$ satisfies \eqref{dissv}, where  $V(t):=\int_0^1 v(t,x)^2 \textrm{dx}$.
\end{corollary}

\subsection{The operators $T$, $K$ an $B$}
Let  $k\in L^2((0,1)^2)$, whose choice will be made precise later on.
 For $f\in L^2(0,1)$, we introduce the transformation 
\begin{equation}\label{T}
g(\cdot)=Tf:x\mapsto f(x)-\int_0^1 k(x,y)f(y) \, \textrm{dy}.
\end{equation}
Notice that for any $f\in L^2(0,1),\, Tf\in L^2(0,1)$. For a solution $u$ of \eqref{system}, we introduce, for any time $t\in [0,\infty)$,
\begin{equation}\label{tf}v(t,\cdot)=Tu(t,\cdot): x\mapsto u(t,x)-\int_0^1 k(x,y)u(t,y) \, \textrm{dy}.\end{equation}
We would like to choose  $k\in L^2((0,1)^2)$ in such a way that  $v$ verifies the target system \eqref{cibleFredholm}, with initial condition $$v_0(x)=u_0(x)-\int_0^1 k(x,y)u_0(y) \, \textrm{dy} \in L^2(0,1).$$
Remark that the condition $v(t,0)=0$ is ensured as soon as $k(0,y)=0$ on $(0,1)$.
We obtain, formally, the equation on the kernel $k$ by taking the time derivative of \eqref{tf}, using \eqref{system} and \eqref{cibleFredholm}, and performing integrations by parts:
\begin{align*}
& (x^\alpha v_x(t,x))_x-\lambda v(t,x)=(x^\alpha u_x(t,x))_x-\int_0^1 k(x,y)(y^\alpha u_y(t,y))_y \, \textrm{dy} \\
\Rightarrow  (x^\alpha u_x(t,x))_x- \lambda u(t,x) -&\int_0^1 \left((x^\alpha k_x(x,y))_x -\lambda k(x,y)\right)u(t,y) \, \textrm{dy}=(x^\alpha u_x(t,x))_x-\int_0^1 k(x,y)(y^\alpha u_y(t,y))_y \, \textrm{dy} \\
\Rightarrow  -&\int_0^1 \left((x^\alpha k_x(x,y))_x -\lambda k(x,y)\right) u(t,y)\, \textrm{dy}+  \lambda u(t,x)=-\int_0^1 (y^\alpha k_y(x,y))_y u(t,y)\, \textrm{dy},
\end{align*}
provided that the kernel $k$ satisfies the boundary conditions $k(0,y)=k(x,0)=k(x,1)=k_y(x,1)=0$. Finally, we impose the kernel $k$ to satisfy the following PDE:
\begin{equation}\label{kernelT1}
\begin{cases}
- (y^\alpha k_y(x,y))_y + (x^\alpha k_x(x,y))_x -\lambda k(x,y) =-\lambda \delta(x=y),  &(x,y)\in (0,1)^2,\\
k(0,y)=0, &y\in (0,1), \\
k(x,0)=k(x,1)=0, &x\in (0,1), \\ 
k_y(x,1)=0,  &x\in (0,1),
\end{cases}
\end{equation}
where $\delta(x=y)$ is the Dirac distribution on the hypersurface $\{(x,y)\in (0,1)^2, x=y\}$. 
Let us decompose formally $k$ (which is assumed to belong to $L^2((0,1)\times(0,1))$) as 
\begin{equation}
\label{decomposition-kernel}
k(x,y)=\sum_{n\in \NN^*} \psi_n(x)\phi_n(y),
\end{equation}
 where, for all $n\in\NN^*$, $\phi_n$ are the eigenfunctions associated with the operator $-A$ introduced in \eqref{ef}.

Therefore \eqref{kernelT1} is equivalent to solve (in appropriate spaces to be defined later on), for any $n\in \NN^*$,
\begin{equation}\label{kernelT2}
\begin{cases}
\lambda_n \psi_n(x) + (x^\alpha \partial_x \psi_n(x))_x -\lambda \psi_n(x) =-\lambda \phi_n(x),  &x\in (0,1),\\
\psi_n(0)=0,  \\ 
\sum_{k\in \NN^*} \psi_k(x)\phi_k'(1)=0, &x\in (0,1).
\end{cases}
\end{equation}
We simplify \eqref{kernelT2} by introducing the change of unknowns
\begin{equation*}
 \psi_n=\phi_n-\xi_n.
\end{equation*}
 Noticing that $\phi_n$ verifies \eqref{ef}, this implies that $\xi_n$ satisfies
\begin{equation}\label{kernelT3}
\begin{cases}
\lambda_n\xi_n(x) + (x^\alpha \partial_x \xi_n(x))_x -\lambda\xi_n(x) =0,  &x\in (0,1),\\
\xi_n(0)=0,  \\ 
\sum_{k\in \NN^*} \xi_k(x)\phi_k'(1)=\sum_{k\in \NN^*} \phi_k(x)\phi_k'(1), &x\in (0,1).
\end{cases}
\end{equation}
The solutions to the first two equations of \eqref{kernelT3} can be written under the form
\begin{equation}\label{pn}
\xi_n(x)=c_n\dfrac{(2\kappa)^{1/2}x^{(1-\alpha)/2}}{J'_{\nu}(j_{\nu,n})}J_{\nu}\left(\dfrac{\sqrt{\lambda_n-\lambda}}{\kappa} x^\kappa\right), \quad x\in (0,1), \quad n\in \NN^*,
\end{equation}
with $c_n\in \RR$ to be determined. Notice that, in the above expression, if $\lambda>\lambda_n$, which might happen only for a finite number of $\lambda_n$, then $\sqrt{\lambda_n-\lambda}$ has to be understood as $i\sqrt{\lambda-\lambda_n}$.
Let us introduce the following function
\begin{equation}\label{wtpsin}
\widetilde{\xi}_n(x)= \dfrac{(2\kappa)^{1/2}x^{(1-\alpha)/2}}{J^\prime_{\nu}(j_{\nu,n})}J_{\nu}\left(\dfrac{\sqrt{\lambda_n-\lambda}}{\kappa} x^\kappa\right), \quad x\in (0,1), \quad n\in \NN^*,
\end{equation}
so that  $\psi_n$ can now be written as 
\begin{equation}
\label{pp} \psi_n=\phi_n-c_n\widetilde{\xi}_n.
\end{equation}
To obtain a solution of  \eqref{kernelT2}, at least formally, it is necessary to deduce the existence of a sequence $c_n$ such that the last equation of \eqref{kernelT3} is satisfied.  In fact, as we will see, our strategy does not require to give an appropriate meaning to the solutions of  \eqref{kernelT1} and to solve explicitly this equation. What is important is to ensure the existence of a solution to \eqref{kernelT2}. We then define $k$ by the formula \eqref{decomposition-kernel}, where $\psi_n$ is given by \eqref{pp} and  $c_n$ is such that \eqref{kernelT3} is verified. The transformation $T$ is then given by \eqref{T}. At this point, we only need to ensure that $k\in L^2((0,1)^2)$. The fact that $Tu$ is indeed a solution of \eqref{cibleFredholm} will be obtained by using an abstract argument based on semigroup theory similar to the one developed in \cite{CGM}.

For our purposes, and following \cite{CGM}, we introduce the following operator:
\begin{equation}\label{K}
K: f\in L^2(0,1) \mapsto \int_0^1 k(1,y)f(y)dy,
\end{equation}
where $k$ is defined in \eqref{decomposition-kernel}. The fact that $K$ is well-defined will be proved later on.
To conclude, let us remind that in this case, an easy duality argument enables us to define our control operator as
\begin{equation}\label{defB}
B: z\in \mathbb R\mapsto-z\delta_1' \in D(A)',
\end{equation}
where $\delta_1$ is the Dirac measure at point $1$. Indeed, recall that \eqref{defB} means that for any $\varphi\in D(A)$, $\varphi'(1)$ exists and for any $z\in \mathbb R$,
$$\langle Bz, \varphi \rangle_{D(A)',D(A)} = z\varphi'(1).$$
With these notations and setting $U=Ku$, equation \eqref{system} can be rewritten  in an abstract way as 
\[
\begin{cases}
\partial_t u=(A+BK)u, \quad t\in (0,T) \\
u(0)=u_0.
\end{cases}
\]

\subsection{Main results}

Following closely what has been done in \cite{CoronLu14}, \cite{CoronLu15}, \cite{CGM}, we need to prove several results in order to establish that our system \eqref{system} is rapidly exponentially stabilized.

\begin{theorem}[\textcolor{black}{Properties of}  \eqref{kernelT3}]
\label{thm-wp}Assume that \eqref{nonresonance} holds.
There exists a unique sequence $(c_n)_{n\in \mathbb N^*}$  such that 
\begin{equation}\label{cndn} c_n-1 \in l^2(\mathbb N^*) \end{equation}
and  such that for any $n\in\mathbb N^*$, the corresponding $\psi_n$ defined in \eqref{pp} verifies \eqref{kernelT2}. Moreover, the corresponding kernel $k$ defined  by the formula \eqref{decomposition-kernel} verifies $k\in L^2((0,1)^2)$.

\end{theorem}

The proof of Theorem \ref{thm-wp} relies on the adaptation of the  finite-dimensional strategy explained and proved in \cite[Section 1.2]{CGM}. It is based on some spectral analysis involving a Riesz basis. In the infinite-dimensional case, the analysis is more complicated, because the control operator needs to be unbounded. It introduces regularity issues. 
\begin{theorem}[Continuity and invertibility of \eqref{T}]
\label{thm-prop}Assume that \eqref{nonresonance} holds. 
For the kernel $k\in L^2((0,1)^2)$ given in Theorem \ref{thm-wp}, the operator $T$ given by \eqref{T} is continuous from $L^2(0,1)$ to $L^2(0,1)$ and invertible.
\end{theorem}

The proof of Theorem \ref{thm-prop} strongly relies on the theoretical operator approach developed in \cite{CGM}. The main difficulty is to prove the invertibility of \eqref{T} by first proving that $T$ is a Fredholm operator and then proving that $T^*$ is one-to-one by studying its kernel.

Thanks to Theorems \ref{thm-wp} and \ref{thm-prop}, we are able to  deduce the rapid exponential stabilization of \eqref{system}. 
\begin{theorem}[Rapid stabilization of \eqref{system}]\label{thm-stab}
For any $\lambda>0$ verifying \eqref{nonresonance}, there exists $C(\lambda)>0$  and a feedback law 
$U(t)=K(u(t)),$
where $K\in L^2(0,1)'$ is given by \eqref{K},
such that for any $u_0\in L^2(0,1)$, there exists a unique solution $u$ of \eqref{system} that verifies moreover: for any $t\geqslant 0$,
$$||u(t,\cdot)||_{L^2(0,1)}\leq C(\lambda)||u_0||_{L^2(0,1)} e^{-\lambda t}.$$
\end{theorem}

Theorem \ref{thm-stab} relies on an abstract strategy based on the semigroup theory together with the invertibility of $T$ and the dissipation estimate \eqref{dissv}. 
\section{Existence of $k$}

\label{sec_wp}

This section is devoted to the proof of Theorem \ref{thm-wp}. It is divided into two subsections: the first one tackles the uniqueness property, while the second one deals with the existence and the desired regularity of $k$. Indeed, the uniqueness will give us instructive informations to treat the existence part of our proof. 

\subsection{Uniqueness}
\label{sec_uniqueness}
Let us first address the uniqueness of the solution to  \eqref{kernelT2}. Using decomposition \eqref{pp}, we consider two sequences $(c_n)_{n\in\mathbb N^*}$ and $(\tilde c_n)_{n\in \mathbb N^*}$ verifying \eqref{cndn}, and such that $\psi_n=\phi_n-c_n\widetilde{\xi}_n$ and $\tilde \psi_n=\phi_n-\tilde c_n\widetilde{\xi}_n$ verify \eqref{kernelT2}. By linearity,
\begin{equation}\label{hatp}\widehat \psi_n:=\psi_n- \tilde \psi_n=\left(\tilde c_n-c_n \right)\widetilde{\xi_n} \end{equation}
verifies
\begin{equation}\label{kernelNul}
\begin{cases}
\lambda_n \widehat \psi_n(x) + (x^\alpha \partial_x \widehat \psi_n(x))_x -\lambda \widehat \psi_n(x) =0,  &x\in (0,1),\\
\widehat \psi_n(0)=0,  \\ 
\sum_{n\in \NN^*} \widehat \psi_n(x)\phi_n'(1)=0,  &x\in (0,1).
\end{cases}
\end{equation}
We can obtain that $\widehat \psi_n \equiv 0$ by proving the following property:
\begin{lemma}\label{Rieszuniq}
The family $\{\sqrt{\lambda_n}\widetilde{\xi}_n\}_{n\in \NN^*}$  is a Riesz basis in $D(A^{\frac{1}{2}})'$.
\end{lemma}
Assume (for the moment) that Lemma \ref{Rieszuniq} holds true. By \eqref{cndn}, we observe that $(\tilde c_n-c_n)\in l^2(\mathbb N^*)$. Since  $\{\sqrt{\lambda_n}\widetilde{\xi}_n\}_{n\in \NN^*}$ is a Riesz basis in $D(A^{\frac{1}{2}})'$, we deduce that $$\sum_{n\in \NN^*}  \left(\tilde c_n-c_n \right)\sqrt{\lambda_n}\widetilde{\xi_n}\in D(A^{\frac{1}{2}})'$$ and 
there exists $C>0$ such that   
\begin{equation}\label{rga}\left |\left |\sum_{n\in \NN^*}  \left(\tilde c_n- c_n \right)\sqrt{\lambda_n}\widetilde{\xi_n}(x)\right | \right |^2_{D(A^\frac{1}{2})'}\geqslant C \sum_{n\in \NN^*} \left(\tilde c_n- c_n \right)^2.\end{equation}
From \eqref{dln} and \eqref{behaviorB}, we have  
 $
\left (\frac{\phi_n'(1)}{\sqrt{\lambda_n}} \right)\in l^\infty(\mathbb{N}^*).$ Hence, writing 
$$\left(\tilde c_n- c_n \right)\widetilde{\xi_n}(x)\phi_n'(1)=\left(\frac{\phi_n'(1)}{\sqrt{\lambda_n}}\right)\left(\tilde c_n- c_n \right) \sqrt{\lambda_n} \widetilde{\xi_n}(x)$$
and remarking that $\left\{\left(\tilde c_n- c_n \right)\frac{\phi_n'(1)}{\sqrt{\lambda_n}}\right\}\in l^2(\mathbb N^*)$,
we deduce by \eqref{hatp} and \eqref{rga} that the last line of 
\eqref{kernelNul} holds if and only if $\left(\tilde c_n- c_n \right)\frac{\phi_n'(1)}{\sqrt{\lambda_n}}=0$ for all $n\in \NN^*$, \textit{i.e.} $c_n=0$ since $\phi_n'(1)\not =0$ by \eqref{behaviorB}. Coming back to \eqref{hatp}, we deduce that $\widehat \psi_n \equiv 0$. Uniqueness of the solution of \eqref{kernelT2} is then proved.
\\

The remaining parts of this section is devoted to the proof of Lemma \ref{Rieszuniq}. This proof relies on two fundamental results on Riesz basis. Before recalling them, let us give a definition.
\begin{definition}[$\omega$-independent sequence]
Let $H$ be a Hilbert space and $\lbrace g_n\rbrace\subset H$. The sequence $\lbrace g_n\rbrace$ is said to be $\omega$-independent if for any sequence $(a_n)_{n\in\mathbb N^*}$ of real numbers,
\begin{equation}
\sum_{n\in\NN^*} a_n g_n = 0 \text{ and }  \sum_{n\in\NN^*}|a_n|^2||g_n||_H^2<\infty \Rightarrow a_n=0,\: \forall n\in\mathbb{N}^*.
\end{equation}
\end{definition}

Let us now state the two abovementioned results on Riesz basis.
\begin{theorem} \cite[Theorem 15]{YoungBook} \label{thm_independent}
Let $H$ be an infinite dimensional separable Hilbert space and let $\{e_n\}_{n\in \NN^*}$ be a Hilbert basis of $H$. If $\{g_n\}_{n\in \NN^*}$ is an $\omega$-independent sequence quadratically close to $\{e_n\}_{n\in \NN^*}$, \textit{i.e.}
$
\sum_{n\in \NN^*}\|e_n-g_n\|^2_{H} < +\infty$, then $\{g_n\}_{n\in \NN^*}$ is a Riesz basis for $H$.
\end{theorem}

\begin{theorem} \label{thm_dense}
Let $H$ be an infinite dimensional separable Hilbert space and let $\{e_n\}_{n\in \NN^*}$ be an Hilbert basis of $H$. If $\{g_n\}_{n\in \NN^*}$ is complete in $H$ and is quadratically close to $\{e_n\}_{n\in \NN^*}$, then $\{g_n\}_{n\in \NN^*}$ is a Riesz basis for $H$.
\end{theorem}
A proof of Theorem \ref{thm_dense}, stated as a remark in \cite[Remark 2.1, p. 318]{Gohberg}, may be found in \cite[Proof of Theorem 3.3]{CGM}. We are now in position to prove Lemma \ref{Rieszuniq}.

\textbf{Proof of Lemma \ref{Rieszuniq}.}
Firstly, we remark that  by the definition of the norm on the space $D(A^{\frac{1}{2}})'$ defined in \eqref{dual12}, $\{\sqrt{\lambda_n} \phi_n\}$ is a Hilbert basis of $D(A^{\frac{1}{2}})'$. Our proof is divided into two steps. In a first step, we will prove that $\{\sqrt{\lambda_n}\widetilde{\xi}_n\}$ is quadratically close to $\{\sqrt{\lambda_n} \phi_n\}$, \textit{i.e.}
\begin{equation}\label{squareint}
\sum_{n\in \NN^*}{\lambda_n}\left\| \phi_n-\widetilde{\xi_n} \right\|_{D\left(A^{\frac{1}{2}}\right)^{\prime}}^2< + \infty. 
\end{equation}
In a second step, we will use some spectral properties of degenerate parabolic equations and apply a contradiction argument in two different cases: the first one will apply Theorem \ref{thm_dense}, while the second one will apply Theorem \ref{thm_independent}. This will allow us to conclude that $\{\widetilde{\xi}_n\}_{n\in \mathbb N^*}$ is a Riesz basis. 
\\

\textbf{First step: $\{\sqrt{\lambda_n}\widetilde{\xi}_n\}$ is quadratically close to $\{\sqrt{\lambda_n}\phi_n\}$.}
We remark that, for a fixed $n\in\mathbb N^*$, we have $\sqrt{\lambda_n}\left( \phi_n-\widetilde{\xi_n} \right)\in L^2(0,1)$, so that it is in $D(A^{\frac{1}{2}})'$. Hence, in order to prove \eqref{squareint}, we can assume that we only consider the $n$ large enough such that $\lambda_n-\lambda>0$, \textit{i.e.} $\sqrt{\lambda_n-\lambda}>0$. Let us prove that 
\begin{equation}\label{squareint2}
\sum_{\lambda_n>\lambda}{\lambda_n}\left\| \phi_n-\widetilde{\xi_n} \right\|_{D\left(A^{\frac{1}{2}}\right)^{\prime}}^2 < + \infty. 
\end{equation}
One has by definition of the norm of $D(A^\frac{1}{2})'$:
\begin{equation}\label{pro}
 \sum_{\lambda_n>\lambda}{\lambda_n}\left\| \phi_n-\widetilde{\xi_n} \right\|_{D\left(A^{\frac{1}{2}}\right)^{\prime}}^2 = \sum_{\lambda_n>\lambda}\sum_{k\in \NN^*} \frac{{\lambda_n}}{\lambda_k}\left| \left< \left( \phi_n-\widetilde{\xi_n} \right), \phi_k \right> \right|^2.
\end{equation}
We introduce $\varepsilon_n$ given by
\begin{equation}\label{defvare}
\varepsilon_n:=j_{\nu,n}-\dfrac{\sqrt{\lambda_n-\lambda}}{\kappa}=j_{\nu,n}-j_{\nu,n}\sqrt{1-\dfrac{\lambda}{(\kappa j_{\nu,n})^2}},
\end{equation}
where we have used the fact that $\lambda_n$ is given by \eqref{dln}. From \eqref{zeroasympt}, we remark that for $n\rightarrow \infty$, $j_{\nu,n}= n + \mathcal{O}(1/n)$. We deduce that 
\begin{equation}\label{den}
\varepsilon_n= \dfrac{\lambda}{2j_{\nu,n}\kappa^2}+\mathcal{O}\left(\dfrac{1}{n^3}\right).
\end{equation}
Notably, by \eqref{zeroasympt}, we deduce that 
\begin{equation}\label{den2}
\varepsilon_n=\mathcal{O}\left(\dfrac{1}{n}\right).
\end{equation}
Using \eqref{IPB}, \eqref{IPB2}, \eqref{PropOrthogonal}, \eqref{dln} and the fact that $j_{\nu,n}$ are zeros of the Bessel function $J_\nu$, we obtain that: 
\begin{equation}
\label{zo}
\begin{aligned}
\left<(\phi_n-\widetilde{\xi_n}), \phi_k \right>=&\:2\kappa  \int_0^1  x^{(1-\alpha)} \left(\dfrac{1}{J'_\nu(j_{\nu,n})} \left[ J_\nu(j_{\nu,n}x^\kappa) \right .\right.\\&\left . \left .- J_\nu\left(\dfrac{\sqrt{\lambda_n-\lambda}}{\kappa}x^\kappa\right)\right] \right) \dfrac{J_\nu(j_{\nu,k}x^\kappa) }{J^\prime_\nu(j_{\nu,k})} \textrm{\textrm{dx}} \\ 
=&\: \delta_{kn} -2\int_0^1 y \left(\dfrac{1}{J^\prime_\nu(j_{\nu,n})}  J_{\nu}\left(\dfrac{\sqrt{\lambda_n-\lambda}}{\kappa}y\right)\right) \dfrac{J_\nu(j_{\nu,k}y) }{J^\prime_\nu(j_{\nu,k})} \textrm{dy}  \\
=  &  \delta_{kn}
-\dfrac{2\kappa^2}{(\lambda_k-\lambda_n+\lambda)J'_\nu(j_{\nu,n})J^\prime_\nu(j_{\nu,k})} \left[\frac{\sqrt{\lambda_n-\lambda}}{\kappa}J_{\nu}(j_{\nu,k})J_{\nu}^\prime\left(\frac{\sqrt{\lambda_n-\lambda}}{\kappa}\right)\right.\\
&\left.-j_{\nu,k}J^\prime_\nu(j_{\nu,k})J_{\nu}\left(\frac{\sqrt{\lambda_n-\lambda}}{\kappa}\right)\right]\\ 
=&   \left( \delta_{kn} +\dfrac{2j_{\nu,k}\kappa^2}{(\lambda_k-\lambda_n+\lambda)J'_\nu(j_{\nu,n})}J_\nu\left(\dfrac{\sqrt{\lambda_n-\lambda}}{\kappa}\right)\right)
\end{aligned}
\end{equation}
Using \eqref{dln} together with the definition of $\varepsilon_n$ given in \eqref{defvare}, we can apply a Taylor expansion and deduce that there exists $e_n\in [j_{\nu,n}-\varepsilon_n,j_{\nu,n}]$ such that 
\begin{equation}\label{plc1}\left(1+ \dfrac{2j_{\nu,n}\kappa^2}{\lambda J'_\nu(j_{\nu,n})}J_\nu\left(\dfrac{\sqrt{\lambda_n-\lambda}}{\kappa}\right)\right)= \left(1+ \dfrac{2j_{\nu,n}\kappa^2}{\lambda J'_\nu(j_{\nu,n})}\left(-\varepsilon_n J'_\nu(j_{\nu,n}) +\frac{\varepsilon_n^2}{2}J''_\nu(e_n) \right)\right).\end{equation}
Using the fact that Bessel functions are solution  of \eqref{edob} and the fact that $j_{\nu,n}$ are zeros of the Bessel functions, we deduce that
\begin{equation}\label{relzero}
j_{\nu,n}^2 J^{\prime\prime}_{\nu}(j_{\nu,n})+j_{\nu,n} J^\prime_{\nu}(j_{\nu,n})=0.
\end{equation}
By the definition of $e_n$ and \eqref{den2}, we obtain
$J_{\nu}^{\prime\prime}(e_n)\sim_{n\rightarrow\infty} J_{\nu}^{\prime\prime}(j_{\nu,n}).$
Hence,  by \eqref{relzero}, \eqref{zeroasympt} and \eqref{dern}, we have 
\begin{equation}\label{acd}
J_{\nu}^{\prime\prime}(e_n)\sim_{n\rightarrow\infty} -\frac{J_{\nu}^{\prime}(j_{\nu,n})}{j_{\nu,n}}.
\end{equation}
Coming back to \eqref{plc1}, using \eqref{acd}, \eqref{den2}, \eqref{zeroasympt} and \eqref{dern}, we deduce that 
\begin{equation}\label{plc2}\left<(\phi_n-\widetilde{\xi_n}), \phi_n \right>=\left(1+ \dfrac{2j_{\nu,n}\kappa^2}{\lambda J'_\nu(j_{\nu,n})}J_\nu\left(\dfrac{\sqrt{\lambda_n-\lambda}}{\kappa}\right)\right)=1-1+\mathcal O\left(\frac{1}{n^2}\right)+\mathcal O\left(\frac{1}{n^2}\right)=\mathcal O\left(\frac{1}{n^2}\right).\end{equation}
Notice also that the above computations easily lead to the following estimation:
\begin{equation}\label{plc22}J_\nu\left(\dfrac{\sqrt{\lambda_n-\lambda}}{\kappa}\right)=\mathcal O\left(\frac{1}{n^{\frac{3}{2}}}\right) \quad\mbox{and} \quad \frac{1}{n^{\frac{3}{2}}}=\mathcal O\left (J_\nu\left(\dfrac{\sqrt{\lambda_n-\lambda}}{\kappa}\right)\right).
\end{equation}
Hence, using  \eqref{plc22}, \eqref{dln}, \eqref{zeroasympt}, \eqref{dern} and \eqref{nonresonance}, we deduce that for $n\not =k$, we have
\begin{equation}\label{plc222}
\begin{aligned}
\left |\left<\frac{\sqrt{\lambda_n}}{\sqrt{\lambda_k}}\left(\phi_n-\widetilde{\xi_n}\right), \phi_k \right> \right |=&\frac{\sqrt{\lambda_n}}{\sqrt{\lambda_k}}\dfrac{2j_{\nu,k}\kappa^2}{(\lambda_k-\lambda_n+\lambda)J'_\nu(j_{\nu,n})} \mathcal{O} \left(\dfrac{1}{n^{\frac{3}{2}}}\right).\\\leqslant &\dfrac{C_1}{|k^2-n^2|},
\end{aligned}\end{equation}
where $C_1$ is a positive constant depending on $\kappa$ and $\lambda$. 
Let us decompose the following sum as follows:
$$\sum_{\substack{
k\in \NN^* \\ n\neq k}}  \frac{1}{|k^2-n^2|^2}=\sum_{1\leqslant k <n}\frac{1}{|k^2-n^2|^2}+\sum_{n< k \leqslant 2n}\frac{1}{|k^2-n^2|^2}+\sum_{k>2n}\frac{1}{|k^2-n^2|^2}.$$
The first sum can be estimated as follows:
$$\begin{aligned}\sum_{1\leqslant k <n}\frac{1}{|k^2-n^2|^2}&=\sum_{1\leqslant k <n}\frac{1}{|k+n|^2\,|k-n|^2}&\leqslant \frac{1}{n^2}\sum_{1\leqslant k \leqslant n-1}\frac{1}{|n-k|^2}&\leqslant  \frac{1}{n^2}\sum_{1\leqslant j \leqslant n-1}\frac{1}{j^2}&\leqslant  \frac{C_3}{n^2},\end{aligned}$$
for some $C_3>0$.
The second sum can be estimated in the same way:
$$\begin{aligned}\sum_{n< k \leqslant 2n}\frac{1}{|k^2-n^2|^2}&=\sum_{n< k \leqslant 2n}\frac{1}{|k+n|^2\,|n-k|^2}&\leqslant \frac{1}{n^2}\sum_{n+1\leqslant k \leqslant 2n}\frac{1}{|n-k|^2}&\leqslant  \frac{1}{n^2}\sum_{1\leqslant j \leqslant n-1}\frac{1}{j^2}&\leqslant  \frac{C_3}{n^2}.\end{aligned}$$
The last sum is easier to estimate and we have
$$\begin{aligned}\sum_{k>2n}\frac{1}{|k^2-n^2|^2}\leqslant\sum_{k>2n}\frac{16}{9k^4}\leqslant  \frac{C_4}{n^3}\end{aligned},$$
for some $C_4>0$.
Combining the previous estimates, we deduce that for some $C_5>0$,
\begin{equation}\label{finsum}\sum_{\substack{n\in \NN^* \\ n\neq k}}  \frac{1}{|k^2-n^2|^2}\leqslant \frac{C_5}{n^2},
\end{equation}
Coming back to \eqref{pro} and using \eqref{plc2} together with \eqref{plc22} and \eqref{finsum}, we obtain that  \eqref{squareint} is verified.
\\

\textbf{Second step: $\{\sqrt{\lambda_n}\widetilde{\xi}_n\}$ is complete or $\omega$-independent in $D(A^{\frac{1}{2}})'$.}

We now follow closely the strategy proposed by \cite{CoronLu14}.
Let us consider some $b$ that satisfies
\begin{equation}\label{leb}
\begin{cases}
(x^\alpha \partial_x b)_x =0,  &x\in (0,1),\\
b(0)=0,\,b(1)=1.
\end{cases}
\end{equation}
One can solve explicitly the above boundary problem and deduce that $b(x)=x^{1-\alpha}$. Notably, $b\in L^2(0,1)$ and we can decompose $b$ as $b=\sum_{n\in \NN^*} b_n \phi_n$, with $(b_n)_{n\in\mathbb N^*}\in l^2(\mathbb N^*)$. The coefficients $b_n$ can be expressed as 
\begin{equation}\label{dbn}b_n=\int_0^1 x^{1-\alpha}\phi_n(x) \, \textrm{\textrm{dx}}=\frac{1}{\lambda_n}\int_0^1 x^{1-\alpha}(x^\alpha \phi_n')'\,\textrm{\textrm{dx}}.\end{equation}
An integration by parts gives
\begin{equation}\label{Abnd}b_n=\frac{1}{\lambda_n}\left([x \phi_n']_0^1-\int_0^1(1-\alpha)\phi_n'\,\textrm{\textrm{dx}}\right)=\frac{\phi_n'(1)}{\lambda_n}-\frac{1-\alpha}{\lambda_n}[\phi_n]_0^1=\frac{\phi_n'(1)}{\lambda_n}.\end{equation}
From \eqref{dern}, we deduce that $\forall n\in\mathbb N^*,\,b_n\not=0$. Now, we remark that by \eqref{wtpsin}, we have 
$$\widetilde{\xi}_n(1)= \dfrac{(2\kappa)^{1/2}}{|J^\prime_{\nu}(j_{\nu,n})|}J_{\nu}\left(\dfrac{\sqrt{\lambda_n-\lambda}}{\kappa}\right), \quad n\in \NN^*.$$
From now on, let us denote 
\begin{equation}\label{betan}
\beta_n:=\dfrac{(2\kappa)^{1/2}}{\left| J^\prime_{\nu}(j_{\nu,n})\right|}J_{\nu}\left(\dfrac{\sqrt{\lambda_n-\lambda}}{\kappa}\right),
\end{equation}
so that $\beta_n=\widetilde{\xi_n}(1)$.
Let us introduce the auxiliary function 
\begin{equation}\label{defeta}\eta_n=\widetilde{\xi_n}-\beta_nb.\end{equation}
From \eqref{kernelT2},\eqref{pp} and \eqref{leb}, it is clear that $\eta_n$ verifies
\begin{equation*}
\begin{cases}
\lambda_n \eta_n(x) + (x^\alpha \partial_x \eta_n(x))_x -\lambda \eta_n(x) =-(\lambda_n-\lambda)\beta_n b,  &x\in (0,1),\\
\eta_n(0)=0,  \\ 
\eta_n(1)=0.
\end{cases}
\end{equation*}
In other words, $\eta_n\in D(A)$ and $\eta_n$ verifies 
\begin{equation}\label{funceta}
-A\eta_n=(\lambda_n-\lambda) \eta_n+(\lambda_n-\lambda)\beta_nb.
\end{equation}
Let us first prove that $\lbrace \widetilde{\xi}_n\}$ is complete or $\omega$-independent in $L^2(0,L)$, which will be enough to prove the desired result.
Let us consider some sequence $(a_n)_{n \in \NN^*} \in \mathbb R^{\mathbb N^*}$ such that 
\begin{equation}\label{l2tx}\sum_{n\in \mathbb N^*} |a_n|^2 ||\tilde \xi_n||_{L^2(0,1)}^2<\infty\end{equation}
and 
\begin{equation}\label{nulwt}
\sum_{n\in \NN^*}a_n \widetilde{\xi}_n=0.
\end{equation}
Using \eqref{wtpsin} together with \eqref{IPB2} together with the change of variables $y=x^\kappa$, we remark that 
$$||\tilde \xi_n||_{L^2(0,1)}^2=\dfrac{1}{|J^\prime_{\nu}(j_{\nu,n})|^2}\left( \left(1-\dfrac{\nu\kappa^2}{\lambda_n-\lambda}\right)\left|J_\nu\left(\dfrac{\sqrt{\lambda_n-\lambda}}{\kappa}\right)\right|^2 + \left|J^\prime_\nu\left(\dfrac{\sqrt{\lambda_n-\lambda}}{\kappa}\right)\right|^2\right).$$
Notably, for $n$ large enough and using \eqref{defvare} and \eqref{den} (which implies that $ J^\prime_\nu\left(\dfrac{\sqrt{\lambda_n-\lambda}}{\kappa}\right) \equiv J^\prime\left(j_{\nu,n} \right)$ as $n\rightarrow \infty$.
$$||\tilde \xi_n||_{L^2(0,1)}^2\geqslant \dfrac{1}{2|J^\prime_{\nu}(j_{\nu,n})|^2}\left|J_\nu'\left(\dfrac{\sqrt{\lambda_n-\lambda}}{\kappa}\right)\right|^2\geqslant \frac{C|J^\prime_{\nu}(j_{\nu,n})|^2}{|J^\prime_{\nu}(j_{\nu,n})|^2}\geqslant C,$$
for some $C>0$ independent on $n$. We deduce by \eqref{l2tx} that  $(a_n)_{n \in \NN^*} \in l^2(\mathbb N^*,\mathbb R)$.  Moreover, from \eqref{betan}, \eqref{plc22}, \eqref{zeroasympt} and  \eqref{dern}, we have 
$\beta_n=\mathcal O\left(\frac{1}{n}\right).$
This means in particular that $(\beta_n)_{n\in\mathbb N^*}\in l^2(\mathbb N^*, \mathbb R)$.
Hence, using \eqref{defeta} and \eqref{nulwt}, we have
\begin{equation}\label{checkdep}
\sum_{n\in \NN^*}a_n {\eta}_n+  \left(\sum_{n\in\mathbb N^*}a_n \beta_n\right)b=0.
\end{equation}
From equation \eqref{funceta}, we deduce that
\begin{equation}\label{goodeta}-A^{-1}\eta_n=\frac{1}{\lambda_n-\lambda}\eta_n-\beta_nA^{-1}b.\end{equation}
Applying $-A^{-1}$ on each side of \eqref{checkdep} and using \eqref{goodeta}, we deduce that
\begin{equation}\label{checkdepbis}
\sum_{n\in \NN^*}\frac{a_n}{\lambda_n-\lambda}{\eta}_n=0.
\end{equation}
Using  \eqref{defeta} in equality  \eqref{checkdepbis}, we deduce that 
\begin{equation}\label{checkdepter}
\sum_{n\in \NN^*}\frac{a_n}{\lambda_n-\lambda} {\widetilde \xi}_n-  \left(\sum_{n\in\mathbb N^*}a_n \frac{\beta_n}{\lambda_n-\lambda}\right)b=0.
\end{equation}
Now, applying $-A^{-1}$ on each side of \eqref{checkdepbis} and using one more time \eqref{goodeta}, 
we deduce that 
\begin{equation}\label{checkdepbis2}
\sum_{n\in \NN^*}\frac{a_n}{(\lambda_n-\lambda)^2} {\eta}_n-  \left(\sum_{n\in\mathbb N^*}a_n \frac{\beta_n}{\lambda_n-\lambda}\right)\left (-A\right)^{-1}b=0.
\end{equation}
Using one more time relation \eqref{defeta} in \eqref{checkdepbis2}, we deduce that 
\begin{equation}\label{checkdepter2}
\sum_{n\in \NN^*}\frac{a_n}{(\lambda_n-\lambda)^2} {\widetilde \xi}_n-\left(\sum_{n\in \NN^*}\frac{a_n\beta_n}{(\lambda_n-\lambda)^2}\right)b-  \left(\sum_{n\in\mathbb N^*}a_n  \frac{\beta_n}{\lambda_n-\lambda}\right)\left (-A\right)^{-1}b=0.
\end{equation}
Therefore, we obtain easily by induction that for any $p\in \NN$, 
\begin{equation}\label{recurrenceuniq}
\sum_{n\in \NN^*}\dfrac{a_n}{(\lambda_n-\lambda)^p} \widetilde{\xi}_n=\sum_{j=1}^p \left(\sum_{n\in \NN^*}\dfrac{a_n\beta_n}{(\lambda_n-\lambda)^j} \right)\left(-A\right)^{j-p}b
\end{equation}
holds (with the usual convention that the sum from $1$ to $0$ equals zero, according to \eqref{nulwt}). 

We distinguish here two cases.
\\

\textbf{First case:} Suppose that 
\[
\sum_{n\in \NN^*}\dfrac{a_n\beta_n}{\lambda_n-\lambda} \neq 0.
\]
Since this coefficient appears in the right-hand side of \eqref{recurrenceuniq}, for any $p\in \NN^*$, for the first term of the sum, we easily deduce by starting from the case $p=1$ and using an induction argument that for any $k\in\mathbb N$, we have 
\begin{equation}\label{A-1subset}
\{A^{-k}b\}_{k\in \NN}\subset \textrm{span}\{\widetilde{\xi}_n\}_{n\in\NN^*}.
\end{equation}
Let us consider some $d\in L^2(0,1)$ such that $d\in \left ({\textrm{span}\{\widetilde{\xi}_n\}_{n\in\NN^*}}\right)^\perp,$ \textit{i.e.}  verifying
\[
\left< h,d\right>_{L^2(0,1)}=0, \quad \forall h \in \textrm{span}\{\widetilde{\xi}_n\}_{n\in\NN^*}.
\]
From \eqref{A-1subset}, we deduce that  $\left<A^{-p} b,d\right>_{L^2(0,1)}=0$, for all $p\in \NN$. Writing $d=\sum_{n\in \NN^*} d_n \phi_n$ for some $(d_n)_{n\in\NN^*}\in l^2(\NN^*)$, we obtain that, for every $p\in \NN$,
\begin{equation}\label{derivativeg}
\sum_{n\in \NN^*} \lambda_n^{-p} b_n d_n=0.
\end{equation}
We define the complex function $G : \CC \rightarrow \CC$ as
\[
z\in \CC \mapsto G(z)=\sum_{n\in \NN^*} d_n b_n e^{z/\lambda_n}\in \CC.
\]
Since $(d_n)_{n\in\mathbb N^*}\in l^2(\mathbb N^*),\,\,(b_n)_{n\in\mathbb N^*}\in l^2(\mathbb N^*)$  and since there exists some positive constant $C$ such that, for any $n\in\mathbb N^*$, we have $\lambda_n\geqslant C$, the function $G$ is an entire function (by uniform convergence on compact sets). Moreover, using \eqref{derivativeg}, we have, for any $p\in \NN$
\[
G^{(p)}(0)=\sum_{n\in \NN^*} \lambda_n^{-p} b_n d_n=0.
\]
Therefore, because $G$ is an entire function, $G\equiv 0$. This exactly means that, for any $n\in\NN^*$, we have $b_n d_n=0$. However, since $b_n\not=0$ for any $n\in\mathbb N^*$, we deduce that $d_n=0$ for any $n\in\NN^*$. Hence, we also obtain, in the Hilbert space $L^2(0,1)$, that
$$\left (\textrm{span}\{\widetilde{\xi}_n\}_{n\in\NN^*}\right)^\perp=\{0\},$$
which exactly means that $\{\widetilde{\xi}_n\}_{n\in\NN^*}$  is complete in $L^2(0,1)$.

Since $L^2(0,1)$  is dense in $D(A^{\frac{1}{2}})'$ and $||\cdot||_{D(A^{\frac{1}{2}})'}^2\leqslant \frac{1}{\lambda_1}|| \cdot||_{L^2(0,1)}^2$, it is easy to deduce that $\{\widetilde{\xi}_n\}_{n\in\NN^*}$  is also complete in the space $D(A^{\frac{1}{2}})'$, so that $\{\sqrt{\lambda_n}\widetilde{\xi}_n\}$ is indeed complete   in $D(A^{\frac{1}{2}})'$ (because both families have the same linear span). Since   $\{\sqrt{\lambda_n}\widetilde{\xi}_n\}$  is quadratically close to $\{\sqrt{\lambda_n}\phi_n\}$ in $D(A^{\frac{1}{2}})'$, it  is a Riesz basis of  $D(A^{\frac{1}{2}})'$  by Theorem \ref{thm_dense}.
\\

\textbf{Second case:} 
Notice that the construction of the first case holds whenever there is a non-zero coefficient in the right-hand side of \eqref{recurrenceuniq}, \textit{i.e}. there exists some $p\in\mathbb N^*$ such that 
\[
\sum_{n\in \NN^*}\dfrac{a_n\beta_n}{(\lambda_n-\lambda)^p}\neq 0.
\]
Hence, we can assume from now on that 
\begin{equation}\label{hatgderiv}
\sum_{n\in \NN^*}\dfrac{a_n\beta_n}{(\lambda_n-\lambda)^p}=0, \quad \forall p \in \NN^*.
\end{equation}
Under the same principle as before, define the complex function $\tilde{G} : \CC \rightarrow \CC$
\[
z\in \CC \mapsto \tilde{G} (z)=\sum_{n\in \NN^*} \dfrac{a_n\beta_n}{\lambda_n-\lambda} e^{z/(\lambda_n-\lambda)}\in \CC.
\]
 Thanks to the definition of $\beta_n$ given in \eqref{betan} together with hypothesis \eqref{nonresonance} and relation \eqref{dln}, we observe that  $\frac{\sqrt{\lambda_n-\lambda}}{\kappa}$ cannot be  a zero of $J_\nu$, so that  $\beta_n\not =0$, for any $n\in\mathbb N^*$.
Moreover, $(\frac{1}{\lambda-\lambda_n})_{n\in\mathbb N^*}\in l^\infty(\mathbb N^*)$ thanks to the non-resonance assumption \eqref{nonresonance} and the fact that $\lambda_n\rightarrow \infty$ as $n\rightarrow \infty$ by \eqref{dln} and \eqref{zeroasympt}. Hence, as before, we can prove easily that the function $\tilde{G} $ is an entire function (by uniform convergence on compact sets). 
Moreover, from \eqref{hatgderiv}, $\tilde{G} ^{(p)}(0)=0, \forall p\in \NN^*$ and therefore $\tilde{G} \equiv 0$. This exactly means that for any $n\in\mathbb N^*$, we have 
$$ \dfrac{a_n\beta_n}{\lambda_n-\lambda}=0.$$

Since $\beta_n\not =0, \forall n \in \NN^*$, we deduce that $a_n=0, \forall n \in \NN^*$, which achieves the proof. Indeed, if \eqref{hatgderiv} is satisfied, we proved that the family $\{\widetilde{\xi}_n\}_{n\in\NN^*}$ is $\omega$-independent in $L^2(0,1)$. We easily deduce that $\{\sqrt{\lambda_n}\widetilde{\xi}_n\}$ is also $\omega$-independent  in $D(A^{\frac{1}{2}})'$. Since   $\{\sqrt{\lambda_n}\widetilde{\xi}_n\}$ is quadratically close to $\{\sqrt{\lambda_n}\phi_n\}$, it is a Riesz basis of  $D(A^{\frac{1}{2}})'$  by Theorem \ref{thm_dense}. 
\cqfd

\subsection{Existence}

To obtain  a solution of  \eqref{kernelT2}, it is necessary to deduce the existence of a sequence $c_n$ such that the last equation of \eqref{kernelT3} is satisfied (remind the definition  \eqref{wtpsin} and the relation \eqref{pp}). As in \cite{CGM}, \cite{CoronLu14} or \cite{CoronLu15}, the difficulty is that the right-hand side of the last equation of \eqref{kernelT3} is not in the appropriate space to use the Riesz basis property. We circumvent this issue by writing the last line of \eqref{kernelT3} as follows:
\begin{equation}\label{TB=Bdecompose}
\sum_{n\in \NN^*} c_n\widetilde{\xi_n}(x)\phi_n^\prime(1)=\sum_{n\in \NN^*} \phi_n^\prime(1) \phi_n(x) = \sum_{n\in \NN^*}\phi_n^\prime(1) \left[\widetilde{\xi_n}(x)+ \phi_n(x)-\widetilde{\xi_n}(x) \right].
\end{equation}
We obtain the decomposition of the sequence $c_n$ by stating the following property:
\begin{equation}\label{Basesproches}
 \sum_{n\in \NN^*}|\phi_n'(1)|^2\left\| \phi_n-\widetilde{\xi_n}\right\|_{D\left(A^{\frac{1}{2}}\right)^\prime}^2 < + \infty,
\end{equation}
which is a consequence of the Riesz basis property \eqref{squareint} together with \eqref{dln} and \eqref{behaviorB} as already explained in Section \ref{sec_uniqueness}.
 Write $c_n=1+d_n$, which reduces \eqref{TB=Bdecompose} to 
\begin{equation}\label{TB=Breste}
\sum_{n\in \NN^*} \phi_n'(1)d_n\widetilde{\xi_n}(x) = \sum_{n\in \NN^*}\phi_n'(1) \left[ \phi_n(x)-\widetilde{\xi_n}(x)\right].
\end{equation}

From \eqref{Basesproches}, \eqref{dln},  \eqref{behaviorB}, the definition of $D(A^{\frac{1}{2}})'$ given in \eqref{dual12} and the Riesz basis provided by Lemma \ref{Rieszuniq}, we know that there exists  $(\tilde d_n)_{n\in \mathbb N^*} \in \ell^2(\mathbb N^*)$ such that 
$$\sum_{n\in \NN^*}\phi_n^\prime(1)\left( \phi_n-\widetilde{\xi_n}\right)=\sum_{n\in \NN^*}\sqrt{\lambda_n}\tilde{d_n}\widetilde{\xi_n},$$
\textit{i.e.}
$$\sum_{n\in \NN^*}\phi_n^\prime (1)\left( \phi_n-\widetilde{\xi_n}\right)=\sum_{n\in \NN^*} \phi_n^\prime(1) \frac{\sqrt{\lambda_n}\tilde {d_n}}{\phi_n^\prime(1)}\widetilde{\xi_n}.$$
Using  \eqref{behaviorB}, \eqref{dln} and \eqref{zeroasympt} , we deduce that $\frac{\sqrt{\lambda_n}}{\phi_n^\prime(1)}=\mathcal O(1)$, so that 
we can set $$d_n= \frac{\sqrt{\lambda_n}\tilde {d_n}}{\phi_n'(1)} \mbox{ and }(d_n)_{n\in \mathbb N^*} \in \ell^2(\mathbb N^*).$$
Hence, in accordance with \eqref{decomposition-kernel} and \eqref{pp}, we define the kernel $k$ as 
\begin{equation}\label{ker2}
k(x,y)=\sum_{n\in \NN^*} \left( \phi_n(x)-\widetilde{\xi_n}(x)- d_n\widetilde   \xi_n(x)\right)\phi_n(y).
\end{equation}
Let us prove that $k\in L^2((0,1)^2)$, which will conclude the proof of Theorem \ref{thm-wp}.  Firstly, we will prove that $\widetilde{\xi}_n$ is quadratically close to $\phi_n$ in $L^2(0,1)^2$, \textit{i.e.}
\begin{equation}\label{squareint3}
\sum_{n\in \NN^*} \| \phi_n - \widetilde{\xi}_n \|^2_{L^2(0,1)} < + \infty. 
\end{equation}

As already explained, we can assume that we only consider the $n$ large enough such that $\lambda_n-\lambda>0$, \textit{i.e.} $\sqrt{\lambda_n-\lambda}>0$. Hence, we will prove that 
\begin{equation}\label{squareint4}
\sum_{\lambda_n>\lambda} \| \phi_n - \widetilde{\xi}_n \|^2_{L^2(0,1)} < + \infty. 
\end{equation}
From \eqref{IPB} and \eqref{IPB2}, the fact $j_{\nu,n}$ is a root of $J_\nu$ and the definition of $\kappa$ given in \eqref{defkappa}, we have 
\begin{align*}
\| \phi_n - \widetilde{\xi}_n \|^2_{L^2(0,1)} =&  \dfrac{2\kappa}{|J'_{\nu}(j_{\nu,n})|^2} \int_0^1 \left| x^{(1-\alpha)/2} \left(J_\nu(j_{\nu,n}x^\kappa)-J_\nu\left(\dfrac{\sqrt{\lambda_n-\lambda}}{\kappa}x^\kappa\right)\right)\right|^2 \textrm{\textrm{dx}} \\
=& \dfrac{2}{|J^\prime_{\nu}(j_{\nu,n})|^2} \int_0^1 y \left|  J_\nu(j_{\nu,n}y)-J_\nu\left(\dfrac{\sqrt{\lambda_n-\lambda}}{\kappa}y\right)\right|^2 \textrm{dy}  \\
=& \dfrac{1}{|J^\prime_{\nu}(j_{\nu,n})|^2} \left(|J^\prime_\nu(j_{\nu,n})|^2 - \dfrac{4j_{\nu,n}\kappa^2}{\lambda}J^\prime_{\nu}(j_{\nu,n})J_{\nu}\left(\dfrac{\sqrt{\lambda_n-\lambda}}{\kappa}\right) \right. \\
& \left. + \left(1-\dfrac{\nu\kappa^2}{\lambda_n-\lambda}\right)\left|J_\nu\left(\dfrac{\sqrt{\lambda_n-\lambda}}{\kappa}\right)\right|^2 + \left|J^\prime_\nu\left(\dfrac{\sqrt{\lambda_n-\lambda}}{\kappa}\right)\right|^2 \right),
\end{align*}
where we used the change of variables $y=x^\kappa$ in the second line.

Using $\varepsilon_n$ introduced in \eqref{defvare} and Taylor expansions, we deduce that there exists $e_n,\,\tilde e_n,\widehat e_n\in [j_{\nu,n}-\varepsilon_n ,j_{\nu,n}]$ such that
\begin{align*}
\| \phi_n - \widetilde{\xi}_n \|^2_{L^2(0,1)} =& \dfrac{1}{|J^\prime_{\nu}(j_{\nu,n})|^2} \left(|J^\prime_\nu(j_{\nu,n})|^2 - \dfrac{4j_{\nu,n}\kappa^2}{\lambda}J^\prime_{\nu}(j_{\nu,n})J_{\nu}\left(j_{\nu,n} -\varepsilon_n\right)\right. \\
& \left. + \left(1-\dfrac{\nu\kappa^2}{\lambda_n-\lambda}\right)\left|J_\nu\left(j_{\nu,n} -\varepsilon_n \right)\right|^2 + \left|J_\nu'\left(j_{\nu,n} -\varepsilon_n\right)\right|^2 \right)
\\=& \dfrac{1}{|J^\prime_{\nu}(j_{\nu,n})|^2} \left(|J^\prime_\nu(j_{\nu,n})|^2 - \dfrac{4j_{\nu,n}\kappa^2}{\lambda}J'_{\nu}(j_{\nu,n})[J^\prime_{\nu}(j_{\nu,n})\varepsilon_n+ J^{\prime\prime}_{\nu}(e_n)\dfrac{\varepsilon_n^2}{2}]\right. \\
& \left. + \left(1-\dfrac{\nu\kappa^2}{\lambda_n-\lambda}\right)\left|J^\prime_\nu(j_{\nu,n})\varepsilon_n+J^{\prime\prime}_{\nu}(\tilde e_n)\dfrac{\varepsilon_n^2}{2}\right|^2 + \left|J^\prime_\nu(j_{\nu,n})+J^{\prime\prime}_\nu(\widehat e_n)\varepsilon_n\right|^2 \right) 
\end{align*}
Using a reasoning already performed to obtain \eqref{acd}, we obtain
\begin{equation}\label{acd2}
J_{\nu}^{\prime\prime}(e_n)=\mathcal O\left(\frac{1}{n^{\frac{3}{2}}}\right),\, J_{\nu}^{\prime\prime}(\tilde e_n)=\mathcal O\left(\frac{1}{n^{\frac{3}{2}}}\right)\mbox{ and } J_{\nu}^{\prime\prime}(\widehat e_n)=\mathcal O\left(\frac{1}{n^{\frac{3}{2}}}\right).
\end{equation}
Using \eqref{zeroasympt}, \eqref{dern}, \eqref{den}, \eqref{den2} and \eqref{acd}  in the previous computations, we deduce that 
\begin{align*}
\| \phi_n - \widetilde{\xi}_n \|^2_{L^2(0,1)} 
=& \dfrac{1}{|J'_{\nu}(j_{\nu,n})|^2} \left(J_\nu'(j_{\nu,n})^2 -2J_\nu'(j_{\nu,n})^2 +\mathcal O\left(\frac{1}{n^3}\right)+ O\left(\frac{1}{n^{\frac{7}{2}}}\right)\right. \\
& \left. + \mathcal O\left(\frac{1}{n^3}\right) +J_\nu'(j_{\nu,n})^2+\mathcal O\left(\frac{1}{n^3}\right)+\mathcal O\left(\frac{1}{n^5}\right)\right)\\
=& \dfrac{1}{|J'_{\nu}(j_{\nu,n})|^2} \mathcal O\left (\frac{1}{n^3} \right ).
\end{align*}
Finally, using once again \eqref{zeroasympt} and \eqref{dern}, we infer that there exists $C>0$ such that for any $n\in \mathbb N^*$,
\begin{equation}\label{auti}
\| \phi_n - \widetilde{\xi}_n \|^2_{L^2(0,1)}  \leqslant \dfrac{C}{n^2},
\end{equation}
which gives \eqref{squareint4} and, hence, proves our result \eqref{squareint3}.
 Since $\{\phi_n\}_{n\in\mathbb{N}}$ is a Hilbert basis of $L^2(0,1)$, we have, by using \eqref{auti},
 $$\begin{aligned}
  ||\phi_n-\widetilde{\xi_n}- d_n\widetilde   \xi_n||^2_{L^2(0,1)}&\leqslant  2 \left ( ||\phi_n-\widetilde{\xi_n}||^2_{L^2(0,1)} +|d_n|^2||\widetilde   \xi_n||_{L^2(0,1)}^2\right)\\& \leqslant  2 \left ( ||\phi_n-\widetilde{\xi_n}||^2_{L^2(0,1)} +2|d_n|^2||\widetilde {\xi_n}- \phi_n||_{L^2(0,1)}^2 +2|d_n|^2||\phi_n||_{L^2(0,1)}^2\right)\\& \leqslant  2 \left (\frac{C}{n^2} +2C\frac{|d_n|^2}{n^2}+2|d_n|^2\right).
  \end{aligned}
 $$
 Hence, we have obtained that there exists $C'>0$ such that for any $n\in \mathbb N^*$,
 \begin{equation}\label{psinte} ||\psi_n||_{L^2(0,L)}^2= ||\phi_n-\widetilde{\xi_n}- d_n\widetilde   \xi_n||^2_{L^2(0,1)}\leqslant C' \left (\frac{1}{n^2}+d_n^2 \right).
 \end{equation}
 From \eqref{ker2} and \eqref{psinte}, we deduce that 
$$||k||_{L^2((0,1)^2)}^2=\sum_{n\in\mathbb N^*} ||\phi_n-\widetilde{\xi_n}- d_n\widetilde   \xi_n||^2_{L^2(0,1)}<\infty.$$
This concludes the proof of Theorem \ref{thm-wp}.
\cqfd

\begin{remark}
Notice that the sequence $\{c_n\}$ does not belong to $ \ell^2(\RR)$ but to $ \ell^\infty(\RR)$. As it will be shown Proposition \ref{Tinvertible}, it is sufficient to prove that the transformation $T$ is continuous from $L^2(0,1)$ into itself. Moreover, if $\{c_n\}$ were to be in $\ell^2(\RR)$, then one would obtain that $T$ is compact and therefore not invertible in $L^2(0,1)$. 
\end{remark}

\section{Properties of the transformation}

\label{sec_prop}

This section is devoted to the analysis and the proof of some properties of the transformation $T$. It is indeed crucial to prove that this transformation $T$ is continuous and invertible to ensure that the original system \eqref{system} is stabilized with the desired decay rate. We start by stating and proving the continuity of the operator and pursue by stating and proving the invertibility of the operator. The invertibility is proved thanks to some special properties satisfied by Fredholm operators, which requires to prove that the transformation $T$ that we consider is Fredholm. 

\subsection{Continuity}

We begin this section by proving that $T$ defined in \eqref{T} verifies $T\in \L(L^2(0,1))$  (the linear continuous functions from $L^2(0,1)$ on itself) and $K$ defined in \eqref{K} is well-defined and verifies  $K\in \L(L^2(0,1);\RR)$. In the following, given a Hilbert space $H$, we denote by $\L(H)$ the space of continous linear operators from $H$ to $H$
\begin{proposition}\label{Tcontinue}
The transformation $T$ belongs to  $\L(L^2(0,1))$.
\end{proposition}

\textbf{Proof of Proposition \ref{Tcontinue}.}
Let $f\in L^2(0,1)$. Then, by \eqref{T}, for $x\in(0,1)$, we have

$$Tf(x)=f(x)-\int_0^1 k(x,y)f(y) \, \textrm{dy}.$$
The identity operator $I_{L^2(0,1)}:L^2(0,1)\rightarrow L^2(0,1)$ is of course continuous from $L^2(0,1)$ into 
$L^2(0,1)$. 
Moreover, since  $||k||_{L^2((0,1)^2)}^2<\infty$, it is standard to deduce that $f\in L^2(0,1)\mapsto \int_0^1 k(\cdot,y)f(y) \, \textrm{dy}$ is a linear continuous operator from $L^2(0,1)$ into $L^2(0,1)$, which concludes the proof.


\cqfd

\begin{lemma}\label{KContinue}
The operator $K$ is well-defined and belongs to $\L(L^2(0,1);\RR)$.
\end{lemma}

\textbf{Proof of Lemma \ref{KContinue}.}
Let $f \in L^2(0,1)$. From \eqref{decomposition-kernel}, we have

\begin{equation}
\label{Kfd}\begin{aligned}
Kf  =\int_0^1 k(1,y)f(y) \textrm{dy}  = \sum_{n\in \NN^*} \psi_n(1) < f , \phi_n >,
\end{aligned}
\end{equation}
as soon as this last quantity turns out to be finite.
Let us prove the following lemma.
\begin{lemma}\label{phi1l2}
The sequence $(\psi_n(1))_{n\in \mathbb N^*}$ belongs to $\ell^2(\NN^*)$.
\end{lemma}
\textbf{Proof of Lemma \ref{phi1l2}.}
Recall that by \eqref{pp}, \eqref{wtpsin} and the definition of $\varepsilon_n$ given in \eqref{defvare}, we have
\begin{align*}
|\psi_n(1)|& =\left|-c_n\dfrac{(2\kappa)^{1/2}}{|J'_{\nu}(j_{\nu,n})|}J_\nu \left(\dfrac{\sqrt{\lambda_n-\lambda}}{\kappa}\right) \right|  \leq \left|-c_n\dfrac{(2\kappa)^{1/2}}{|J'_{\nu}(j_{\nu,n})|}J'_{\nu}(j_{\nu,n})\varepsilon_n \right| 
\leq C |c_n| \varepsilon_n.
\end{align*}
Since the sequence $(\varepsilon_n)_{n\in\mathbb N^*}$ belongs to $ \ell^2(\mathbb N^*)$ by \eqref{den2} and $(c_n)_{n\in\mathbb N^*}$ belongs to $l^\infty(\mathbb N^*)$ by \eqref{cndn}, this concludes our proof of Lemma \ref{phi1l2}.

\cqfd 
From Lemma \ref{phi1l2},  $\{\psi_n(1)\}_{n\in\mathbb N^*} \in \ell^2(\RR)$ and  we also have $(< f , \phi_n >)_{n\in\mathbb N^*} \in \ell^2(\RR)$ since $f\in L^2(0,1)$. Therefore, we deduce  by the Cauchy-Schwarz inequality that $K$ belongs to $\L(L^2(0,1);\RR)$. This concludes the proof of Lemma \ref{KContinue}.
\cqfd


\subsection{Operator equalities}

In this section we provide the functional framework for the equalities
$$T(A+BK)=(A-\lambda I)T \,\, \mbox{ and } \,\,
TB=B.$$
We recall that these equalities will be used later on to prove the invertibility.
 
\begin{proposition}\label{ptbb}
The operator $T$ defined by \eqref{tf} and \eqref{decomposition-kernel} can be uniquely extended as a linear continuous operator from $D(A)'$ to $D(A')$ verifying the functional identity $TB=B$ in $D(A)'.$
\end{proposition}
\textbf{Proof of Proposition \ref{ptbb}.}
For $f\in L^2(0,1)=\sum_{n\in\mathbb N^*}a_n \phi_n$, by \eqref{T} and \eqref{decomposition-kernel}, we have 
$$Tf=f -\sum_{n\in \mathbb N^*}a_n \psi_n.$$
Hence, for $f=\sum_{n\in\mathbb N^*}a_n \phi_n\in D(A)'$, it is reasonable to define $Tf$ as
$$Tf=f -\sum_{n\in \mathbb N^*} a_n\psi_n.$$
We have:
\begin{equation}
\label{D(A)'}
\left\Vert \sum_{n\in \mathbb N^*}a_n\psi_n\right\Vert_{D(A)'}^2= \sum_{k\in \mathbb N^*} \frac{1}{\lambda_k^2} \left (\sum_{n\in \mathbb N^*}\langle a_n \phi_k,\psi_n\rangle_{L^2(0,1)}  \right)^2.
\end{equation}
Recalling that $\psi_n = \phi_n-(1+d_n)\widetilde \xi_n$, and following the same strategy than the one used to obtain \eqref{zo}, we end up with 
\begin{equation}\label{ttze}
\langle \phi_k,\psi_n\rangle_{L^2(0,1)} = \delta_{kn} + (1+d_n)\frac{2j_{\nu,\kappa}}{(\lambda_k-\lambda_n+\lambda)J'_{\nu}(j_{\nu,n})} J_\nu\left(\frac{\sqrt{\lambda_n-\lambda}}{\kappa}\right)
\end{equation}
Concerning the case $n=k$, by \eqref{plc2} and the fact that $\langle \phi_n, \phi_n \rangle=1$, we deduce notably that 
$$\langle \phi_n,\tilde \xi_n\rangle_{L^2(0,1)}=1+\mathcal O \left (\frac{1}{n^2} \right )=\mathcal O(1).$$
Combining this with \eqref{plc2}, we deduce that 
$$
\langle \phi_n,\psi_n\rangle_{L^2(0,1)} =\mathcal O\left(\frac{1}{n^2}\right)+ d_n \mathcal O(1)=\mathcal O(1),$$
since $l^2(\mathbb N^*)\subset l^\infty(\mathbb N^*)$.
By the above estimation, we have for some $C_0>0$,
$$\begin{aligned}\sum_{n\in \mathbb N^*} \frac{1}{\lambda_n^2} a_n^2\langle \phi_n,\psi_n\rangle_{L^2(0,1)}^2&\leqslant  
C_0\sum_{n\in \mathbb N^*}  \frac{a_n^2}{n^4}= C_0||f||^2_{D(A)'}<\infty\end{aligned}$$
by the definition of $D(A)'$ given in \eqref{dap}.
\\

Let us treat now the case where $n\neq k$. In the sequel, we denote by $C_1,\ldots$ various constants not depending on $n$ nor $k$. From \eqref{plc222} and \eqref{ttze}, one has
\begin{equation}
\langle \phi_k,\psi_n\rangle_{L^2(0,1)} \leq C_1 (1+d_n) \frac{k}{n(|k^2-n^2|)} 
\end{equation}

Then, one has using the Cauchy-Schwarz inequality
\begin{equation}\label{neq}
\begin{split}
\sum_{k\in \mathbb N^*} \frac{1}{\lambda_k^2} \left (\sum_{n\neq k}a_n\langle \phi_k,\psi_n\rangle_{L^2(0,1)} \right)^2 \leq &\sum_{k\in \mathbb N^*} \frac{1}{\lambda_k^2} \left (\sum_{n\neq k}  \frac{a_n}{n^2}(1+d_n)C_1\frac{k n}{|k^2-n^2|} \right)^2\\
\leq & \left (\sum_{k\in \mathbb N^*} \frac{1}{\lambda_k^2} C_1^2 k^2 \sum_{n\neq k}\frac{n^2}{|k^2-n^2|^2} \right)\left(\sum_{n\neq k} \frac{a_n^2}{n^4}\right).
\\\leq & C_2 \left (\sum_{k\in \mathbb N^*} \frac{1}{\lambda_k^2}  k^2 \sum_{n\neq k}\frac{n^2}{|k^2-n^2|^2} \right)||f||_{D(A)'}^2,
\end{split}
\end{equation} 
since $(1+d_n)_{n\in\mathbb{N}}\in l^\infty(\mathbb{N}^*)$ and by definition of the norm $||\cdot||_{D(A)'}$ given in \eqref{dap}.

It remains to estimate
$$\sum_{k\in \mathbb N^*} \frac{1}{\lambda_k^2} C_1^2 k^2 \sum_{n\neq k}\frac{n^2}{|k^2-n^2|^2}.$$
We remark that 
\begin{equation}
\sum_{\underset{n \neq k}{n\in\mathbb{N}^*}} \frac{n^2}{|k^2-n^2|^2} = \sum_{1<n< k} \frac{n^2}{|k^2-n^2|^2} + \sum_{k< n\leq 2k} \frac{n^2}{|k^2-n^2|^2} + \sum_{n>2k} \frac{n^2}{|k^2-n^2|^2}.
\end{equation}
We have
\begin{align*}
\sum_{1\leq n< k} \frac{n^2}{|k^2-n^2|^2} = \sum_{1\leq n< k} \frac{n^2}{|k-n|^2|k+n|^2}
\leq  \sum_{1\leq n<k}\frac{1}{|k-n|^2}
\leq \sum_{1\leq j<k-1}\frac{1}{j^2}
\leq  {C_3}
\end{align*}
The same reasoning leads to 
\begin{align*}
\sum_{k< n< 2k} \frac{n^2}{|k^2-n^2|^2} \leq {C_3}.
\end{align*}
The last term verifies
\begin{align*}
\sum_{n>2k} \frac{n^2}{|k^2-n^2|^2} 
\leq \sum_{n>2k} \frac{16n^2}{9n^4}
\leq \frac{C_4}{k}.
\end{align*}
Finally, gathering all the bounds given above, one ends up with
\begin{equation*}
\sum_{\underset{n \neq k}{n\in\mathbb{N}^*}} \frac{n^2}{|k^2-n^2|^2} \leq C_5.
\end{equation*}
Recalling \eqref{dln} and \eqref{zeroasympt}, we can deduce that when $n\neq k$, 
\begin{equation}\label{ndk}
\sum_{k\in \mathbb N^*} \frac{1}{\lambda_k^2} \left (\sum_{n\in \mathbb N^*}a_n\langle \phi_k,\psi_n\rangle_{L^2(0,1)} \right)^2\leqslant  C_6\sum_{k\in \mathbb N^*}\frac{1}{k^2}||f||_{D(A)'}\leqslant C_7||f||_{D(A)'}.\end{equation}
Combining \eqref{neq} and \eqref{ndk},  we deduce that  $T:D(A)'\rightarrow D(A)'$ and $T$ is continuous. The fact that the extension is unique is now straightforward by density of $L^2(0,1)$ in $D(A)'$. We now turn to the equality $TB=B$ in $D(A)'$. First notice that 
\[
T^*f=f-\int_0^1 k(y,x)f(y) dy,
\]
Therefore, for $B\in D(A)'$ defined by \eqref{defB} and using the representation of the kernel \eqref{decomposition-kernel}, we have, for any $n\in\mathbb N^*$,
$$\langle TB, \phi_n \rangle_{D(A)',D(A)}= \phi_n'(1)-\sum_{n\in \mathbb N^*}\langle \psi_n,\phi_n \rangle_{L^2(0,1)} \phi_k '(1).$$
From the last line of \eqref{kernelT2}, we deduce that
$$\langle TB, \phi_n \rangle_{D(A)',D(A)}= \phi_n'(1)=\langle B, \phi_n \rangle_{D(A)',D(A)},$$
from which we deduce our result.
\cqfd

 Let us now introduce the space $$H^1_{\alpha,L}(0,1):=\{f\in H^1_{\alpha}(0,1) \, | \, f(0)=0\}$$ and 
 $$D(A)_{L}:=\{f\in H^1_{\alpha,L}(0,1) \, | \, x^\alpha f_x \in H^1(0,1) \}.$$
Finally, we define 
\[
D(A+BK)=\{f\in D(A)_L \, | \, f(1)=Kf\},
\]
where $K$ has been introduced in \eqref{K}.
Let us recall that, for $f\in D(A+BK)$ and $g\in L^2(0,1)$, $(A+BK)f=g$ is equivalent to 
\begin{equation}\label{A+BK}
\begin{cases}
\left(x^\alpha \partial_x f\right)_x =g, \quad x\in (0,1) \\
f(0)=0, \\
f(1)=K(f).
\end{cases}
\end{equation}

We now turn to the operator equality. 

\begin{proposition}\label{opidentite}
For $f\in D(A+BK)$, we have 
\[
T(A+BK)f=(A-\lambda I)Tf \textrm{ in } L^2(0,1).
\]
\end{proposition}

\textbf{Proof of Proposition \ref{opidentite}.}

\color{black}

Let $f\in D(A+BK)$. We recall that $(A+BK)f=g$, for $g\in L^2(0,1)$ if and only if \eqref{A+BK} is satisfied. Note first that, using the last line of \eqref{kernelT2}, which lies in $D(A)'$, and making the duality product with $\phi_k$ where $k\in\mathbb{N}^*$, one can deduce that
\begin{equation}
\label{TB=Bgreat}
\sum_{n\in \NN*} \<\psi_n,\phi_k\> \phi_n'(1) =0, \quad \forall k\in \NN^*
\end{equation}

Moreover, by integrations by parts and using the last line of \eqref{A+BK}, we obtain that
\begin{equation}
\label{relphiA+BK}
\<(A+BK)f, \phi_n\>=\<(x^\alpha \partial_x f)_x,\phi_n\> =\<f,(x^\alpha \partial_x \phi_n)_x\>-K(f)\phi_n'(1)=-\lambda_n \<f,\phi_n\>-K(f)\phi_n'(1).
\end{equation}
Therefore, the sequence $(\lambda_n \<f,\phi_n\>+K(f)\phi_n'(1))_{n\in \NN^*}$ belongs to $\ell^2(\NN^*)$. Moreover, using \eqref{T} and \eqref{decomposition-kernel}, we have 
\begin{align}
\langle T(A+BK)f,\phi_k\rangle &=\langle (A+BK)f,\phi_k\rangle-\sum_{n\in \NN} \langle \psi_n,\phi_k\rangle \<(A+BK)f, \phi_n\> \nonumber \\
&=\sum_{n\in \NN} \<(A+BK)f, \phi_n\> \langle \phi_n,\phi_k\rangle -\sum_{n\in \NN} \psi_n \<(A+BK)f, \phi_n\> \nonumber \\
&= -\lambda_k \langle f,\phi_k\rangle - K(f) \phi_k'(1) +\sum_{n\in \NN}\langle \psi_n,\phi_k\rangle \lambda_n \<f,\phi_n\> \label{equTA+BK},
\end{align}
the last line being due to \eqref{TB=Bgreat} and \eqref{relphiA+BK}.

Notice that for $f\in D(A+BK)$, we have $Tf \in D(A)$ from \eqref{kernelT2} and the definition of $K$. Moreover, using firstly the definition of $A$ and performing an integration by parts, and secondly equation \eqref{kernelT2}, we obtain the two following identities
\begin{equation*}
\begin{split}
\<A\psi_n,\phi_k\>=&-\lambda_k \<\psi_n,\phi_k\> - \psi_n(1)\phi_k'(1), \\
\<A\psi_n,\phi_k\>=& (\lambda -\lambda_n )\<\psi_n,\phi_k\>  -\lambda \delta_{nk}.
\end{split}
\end{equation*}
This implies that
\begin{align*}
\<(A-\lambda I)Tf ,\phi_k\>=&-\<Tf ,(\lambda_k+\lambda) \phi_k\> \\
=&- (\lambda_k+\lambda) \<f,\phi_k\> + \sum_{n\in \NN^*} (\lambda_k+\lambda)\<\psi_n,\phi_k\>\<f,\phi_n\>  \\
=& - (\lambda_k+\lambda) \<f,\phi_k\> + \sum_{n\in \NN^*} (\lambda_n\<\psi_n,\phi_k\> -  \psi_n(1)\phi_k'(1) +  \lambda \delta_{nk})\<f,\phi_n\> \\
=& -\lambda_k \<f,\phi_k\> - K(f)\phi_k'(1) + \sum_{n\in \NN^*} \lambda_n\<\psi_n,\phi_k\>\<f,\phi_n\> 
\end{align*}
From this equation and \eqref{equTA+BK}, one can deduce that, for all $\phi_k\in D(A)$ and all $f\in D(A+BK)$, $\langle (A-\lambda I)Tf,\phi_k\rangle = \langle T(A+BK)f,\phi_k\rangle$, which concludes the proof of the desired result. 
\color{black}
\cqfd


\subsection{Invertibility of the transformation}

To prove that the transformation $T$ is invertible, we proceed in several steps: firstly, we prove that $T$ is Fredholm and, secondly, we use the Fredholm property of $T$ to prove that $\mathrm{Ker}\:T^*$ is reduced to $0$, where $T^*$ denotes the adjoint operator of $T$. This implies the invertibility of the operator by well-known functional analysis arguments.

The first lemma that we state is equivalent to say that the operator $T$ is Fredholm, \textit{i.e.} that $T$ is composed by the sum of a compact operator and an invertible operator.  

\begin{lemma}\label{TFredholm}
The operator $T$ is written as $T=\tilde{T}+C$ where $\tilde{T}$ (resp. $C$) is an invertible (resp. compact) operator from $L^2(0,1)$ to $L^2(0,1)$. Consequently, $T$ is a Fredholm operator of order $0$. 
\end{lemma}

\textbf{Proof of Lemma \ref{TFredholm}.}

 Recall that $k$ can be decomposed as 
\begin{equation}\label{goodk}k(x,y)=\sum_{k=1}^\infty \left( \phi_n(x)-c_n\widetilde{\xi}_n(x)\right)\phi_n(y),\end{equation}
where $c_n=1+d_n$ for some $\{d_n\}_{n\in \NN^*} \subset \ell^2(\RR)$.
Let $f\in L^2(0,1)$. From the definition of $T$, 
\begin{align*}
Tf&= f-\int_0^1 k(x,y)f(y)\textrm{dy}\\
& =  f- \sum_{n\in \NN^*} (\phi_n(x)-c_n\widetilde{\xi}_n(x))<f,\phi_n> \\
& =  \sum_{n\in \NN^*} (1+d_n) \widetilde{\xi}_n<f,\phi_n>,
\end{align*}
We define 
\[
\tilde{T}f:=\sum_{n\in \NN^*} \widetilde \xi_n(x)<f,\phi_n>, \quad Cf:=\sum_{n\in \NN^*} d_n\widetilde \xi_n(x)<f,\phi_n>.
\]
The continuity from $L^2(0,1)$ to itself  of $\tilde T$ and $C$ follows from the proof of the continuity of $T$ (Lemma \ref{Tcontinue}). 

$\tilde{T}$ is clearly invertible since $\{ \widetilde \xi_n \}_{n\in \NN^*} $ is a Riesz basis of $L^2(0,1)$, Indeed, there exists a constant $C_1>0$ such that for any $f\in L^2(0,1)$, 
\[
\|\tilde{T}f\|_{L^2(0,1)}^2=\left\| \sum_{n\in \NN^*} \widetilde \xi_n(x)<f,\phi_n> \right\|_{L^2(0,1)}^2 \geq C_1 \sum_{n\in \NN^*} |<f,\phi_n>|^2 \geq C_1 \|f\|_{L^2(0,1)}^2,
\]
so that $\tilde{T}$ is injective, and the surjectivity can be deduced by the fact that $\{ \widetilde \xi_n \}_{n\in \NN^*} $ is a Hilbert basis.

We prove the compactness of $C$ by proving that $C$ is a Hilbert-Schmidt operator, that is, 
\[
\sum_{n\in \NN^*} \left\| C \phi_n \right\|_{L^2(0,1)}^2 < +\infty. 
\]
Indeed, 
\[
\sum_{n\in \NN^*} \left\| C \phi_n \right\|_{L^2(0,1)}^2 = \sum_{n\in \NN^*} \left\| d_n \widetilde \xi_n \right\|_{L^2(0,1)}^2 \leq C_2 \sum_{n\in \NN^*} | d_n |^2 < + \infty.
\]
Then, $\tilde{T}$ is a Fredholm operator of order $0$ since it is invertible. It is well known that $\tilde{T} + C$ remains a Fredholm operator of order $0$, which concludes our desired result. 

\cqfd

We are now in position to prove the invertibility of $T$. Indeed, as mentioned before, proving the invertibility of a Fredholm operator reduces to studying the kernel of its adjoint operator. Let us state our invertibility result:
\begin{proposition}
\label{Tinvertible}
The transformation $T$ is invertible from $L^2(0,1)$ to $L^2(0,1)$. 
\end{proposition}

\textbf{Proof of Proposition \ref{Tinvertible}.}
Since $T$ is a Fredholm operator of order $0$ (Lemma \ref{TFredholm}), we reduce the proof of the invertibility of $T$ to proving $\textrm{Ker }T^* = \{0\}$, where we recall that $T^*$ denotes the adjoint operator of $T$. 

We use the operator equality in a slightly different form 
\[
T(A+BK+(\lambda + \rho) I) = (A+\rho I)T,
\]
where $\rho \in \CC$ is given by the following lemma.
\begin{lemma}\label{rhoinv}
There exists $\rho  \in \CC$ such that $(A+BK+(\lambda + \rho) I)$ is an invertible operator from $D(A+BK)$ to $L^2(0,1)$ and $(A+\rho I)$ is an invertible operator from $D(A)$ to $L^2(0,1)$. 
\end{lemma}

We postpone the proof of Lemma \ref{rhoinv} to the appendix as it is only technical and follows the proof of \cite[Proposition 4.6]{CGM}. From Lemma \ref{rhoinv}, we can write the operator equality under the form 
\[
(A+\rho I)^{-1}T=T(A+BK+(\lambda + \rho) I)^{-1}.
\]
Let $\chi \in \textrm{Ker }T^*$. Then, for all $\phi \in L^2(0,1)$, 
\begin{align*}
0&=\<(A+\rho I)^{-1}T\phi-T(A+BK+(\lambda + \rho) I)^{-1}\phi,\chi\>_{L^2} \\
&=\<\phi,T^*((A+\rho I)^{-1})^* \chi\>-\<(A+BK+(\lambda + \rho) I)^{-1}\phi,T^* \chi\>_{L^2} \\
&=\<\phi,T^*((A+\rho I)^{-1})^* \chi\>_{L^2}.
\end{align*}
Since $\textrm{Ker }T^*$ is finite-dimensional, then it is stable by $((A+\rho I)^{-1})^*$. Therefore, if $\textrm{Ker }T^* \neq \{0\}$, then $((A+\rho I)^{-1})^*$ has an eigenfunction in $\textrm{Ker }T^*$. Moreover, this eigenfunction is also an eigenfunction of $(A^*)^{-1}=A^{-1}$. Therefore this eigenfunction can be chosen on the form $\phi_k$ for a given $k\in \NN$. From the $TB=B$ condition, we deduce 
$
B^* \phi_k= B^* T^* \phi_k = 0$,
which is a contradiction with the fact that
$B^* \phi_k= \phi_k'(1) \neq 0,\,\forall k\in \NN^*$
by \eqref{dern}. Hence, $\mathrm{Ker}\:T^* =\lbrace 0\rbrace$, which implies that $T$ is invertible. This concludes the proof of Proposition \ref{Tinvertible}.

\cqfd

Notice that the proof of Theorem \ref{thm-prop} follows from Propositions \ref{Tcontinue} and \ref{Tinvertible}.

\section{Proof of Theorem \ref{thm-stab}}
\label{ABKS}

\textbf{Proof of Theorem \ref{thm-stab}}


To prove Theorem \ref{thm-stab}, we aim at applying classical results in the semigroup theory. The idea is to prove that the operator $A+BK$ is a $m$-dissipative operator in an appropriate equivalent norm to $||\cdot||_{L^2(0,1)}$. We can deduce then that there exists a unique solution to \eqref{system}. Then, using the invertibility and the continuity of the operator $T$, we can deduce the exponential stability estimate given in Theorem \ref{thm-stab}. We firstly prove that the operator $A+BK$ is dissipative in an appropriate equivalent norm to $||\cdot||_{L^2(0,1)}$ and secondly that it is maximal. Third, we finally prove the exponential stability estimate given in Theorem \ref{thm-stab}.\\

\textbf{First step: $A+BK$ is dissipative.}

Since the operator \eqref{T} given by Theorem \ref{thm-prop} is continuous, $\Vert \cdot\Vert_{T}:=\Vert T \cdot\Vert_{L^2(0,1)}$ is equivalent to the usual norm $\Vert \cdot \Vert_{L^2(0,1)}$. We denote by $\langle \cdot,\cdot\rangle_T$ the associated inner product. Then, proving that $A+BK$ reduces to proving that, for all $u\in D(A+BK)$, $\langle (A+BK)u,u\rangle_{T}\leq 0$. Indeed, we have
\begin{equation}
\begin{split}
\langle T(A+BK)u,Tu\rangle_{L^2(0,1)} = &\langle (A-\lambda I_{L^2(0,1)})Tu,Tu\rangle_{L^2(0,1)} \\
=& \langle ATu,Tu\rangle_{L^2(0,1)} - \lambda \Vert u\Vert^2_{T}\\
\leq &\hspace{0.2cm} 0,
\end{split}
\end{equation}
where, in the first line, we have used Proposition \ref{opidentite} and, in the last line, we have used the fact the operator $A$ is dissipative. \\

\textbf{Second step: $A+BK$ is maximal.} 

Note that the inclusion $\mathrm{Ran}(\sigma I_{L^2(0,1)} - (A+BK))\subset L^2(0,1)$ for some $\sigma>0$ is obvious. We aim at proving that, for some positive constant $\sigma$, $L^2(0,1)\subset \mathrm{Ran}(\sigma I_{L^2(0,1)}  - (A+BK))$. It is equivalent to prove that, for each $u\in L^2(0,1)$, there exists $\tilde{u}\in D(A+BK)$ such that
\begin{equation}
(\sigma I_{L^2(0,1)}-(A+BK)) \tilde{u} = u.
\end{equation}
This equation is obviously equivalent to $T^{-1}T(A+BK-\sigma I_{L^2(0,1)})\tilde{u} = u$. Using Proposition \ref{opidentite}, we thus have
\begin{equation}
(\sigma I_{L^2(0,1)} -(A-\lambda I_{L^2(0,1)})) T \tilde{u}= Tu.  
\end{equation}
Hence, we have proved that the operator $A-\lambda I_{L^2(0,1)}$ is a $m$-dissipative operator, so that it is the generator of a semigroup of contractions by \cite[Proposition 3.1.13]{TucsnakWeissBook} (remark that notably, we automatically have that $D(A+BK)$ is dense in $L^2(0,1)$ by \cite[Proposition 3.1.6]{TucsnakWeissBook}). In particular, \cite[Proposition 2.3.5]{TucsnakWeissBook} holds and we deduce the existence and uniqueness of a solution to \eqref{system} with $U(t)=Ku(t,.).$ 
\\

\textbf{Third step: Exponential stability.} 

Using the continuity and the invertibility of the operator \eqref{T} given by Theorem \ref{thm-prop}, the relation \eqref{tf}, and Corollary \ref{cdv}, one immediately obtains that, for all $t\geq 0$,
\begin{equation}
\begin{split}
||u(t,\cdot)||_{L^2(0,1)} = &\Vert T^{-1} T u(t,\cdot)\Vert_{L^2(0,1)}\\
\leq & \Vert T^{-1}\Vert_{\L(L^2(0,1))} \Vert Tu(t,\cdot)\Vert_{L^2(0,1)}\\ 
\leq & \Vert T^{-1}\Vert_{\L(L^2(0,1))} e^{-\lambda t} \Vert T u_0\Vert_{L^2(0,1)}\\
\leq & \Vert T^{-1}\Vert_{\L(L^2(0,1))}\Vert T\Vert_{\L(L^2(0,1))} e^{-\lambda t} \Vert u_0\Vert_{L^2(0,1)}
\end{split}
\end{equation}
This concludes the proof of the exponential stability given in Theorem \ref{thm-stab}.

\color{black}

\cqfd

\section{Conclusion and open problems}

In this paper, we have solved the rapid stabilization problem of a degenerate parabolic equation. To do so, we have applied a backstepping strategy with a Fredholm transform, as it has been done for instance for the Schr\"odinger equation \cite{CGM}. To do so, we have proved existence and uniqueness of the involved kernel by applying some results about Riesz bases and applied operator theory to deduce properties on the Fredholm transform.\\

This paper paves the way to many other open problems. It is well-known in the literature that the strong degenerate case, \textit{i.e.} $\alpha\geq 1$, is way more difficult, and it could be interesting to address the rapid stabilization problem for this equation. Another issue to tackle is the case where the control is localized at $x=0$. In this case, different Bessel functions should be studied. They are more difficult to handle since they explode at $x=0$. It would also be interesting to treat more general degenerate problems in the spirit of \cite{MR3456387}. Finally, solving the finite-time stabilization problem as in \cite{coron2017null} is also an interesting research line, since it will allow us to solve the null-controllability problem for parabolic degenerate equations. It is a challenging problem, since it is not known nowadays whether a finite-time stabilization can be achieved with parabolic equations and Fredholm transforms.\\

More generally, the use of the Fredholm transform in the backstepping approach seems to exhibit common features when applied to different PDEs. This comes from the fact that one can exploit the properties of the eigenfunction of the spatial differential operator to study the different properties of the kernel. This is not a feature shared for the Volterra transformation for instance. An interesting open problem therefore lies in seeking necessary or sufficient conditions for the existence, continuity and invertibility of a Fredholm transform for the backstepping approach in an abstract setting.


\appendix

\section{Proof of Lemma \ref{rhoinv}}

We now turn to the proof of Lemma \ref{rhoinv}, which follows the proof of \cite[Proposition 4.6]{CGM}. 
Let $\kappa:=\lambda + \rho$. Applying $A^{-1}$ to $A+BK+\kappa I$ yields 
\[
I+A^{-1}BK+\kappa A^{-1} : D(A+BK) \rightarrow D(A).
\] 
Firstly, we are going to  prove that 
$$\{\kappa\in\mathbb C,\quad I+A^{-1}BK+\kappa A^{-1}\mbox{  is invertible}\}$$
is nonempty. Let us first deal with the case where $K(A^{-1}B)\neq -1$. In this case, the operator $I+A^{-1}BK$ is invertible. Indeed, to solve
\[
(I+A^{-1}BK)f=g,
\]
for any $g\in D(A)$, one just has to apply $K$ to the last equation, which yields 
\[
K(f)(1+K(A^{-1}B))=K(g),
\]
which implies, using $1+K(A^{-1}B)\neq 0$, 
\[
f=g-\dfrac{A^{-1}BK(g)}{1+K(A^{-1}B)}.
\]

Assume now that $K(A^{-1}B)=-1$. In this case, since $A^{-1}B(x)= x^{1-\alpha}$ by \eqref{leb}, we have $A^{-1}B \in D(A+BK)$. Notice that $0$ is an eigenvalue of $I+A^{-1}BK$ associated to the eigenvector $A^{-1}B$, of algebraic multiplicity $1$ ($f+A^{-1}BKf=0$ implies that $f$ is colinear to $A^{-1}B$). Then, from \cite{Nagy}, there exists an open set $\Omega \subset \CC $ of $0\in \CC$ such that there exist an holomorphic function $\kappa \in \Omega \mapsto \lambda(\kappa) \in \CC$ and an holomorphic function $\kappa \in \Omega \mapsto f(\kappa) \in D(A+BK)$ such that
\begin{align}
f(0)&=A^{-1}B, \qquad \nonumber \\
(I+A^{-1}BK+\kappa A^{-1})f(\kappa)&=\lambda(\kappa) f(\kappa),\quad \kappa\in\Omega. \label{vpinv}
\end{align}
Let us use a proof by contradiction and assume that $\lambda(\kappa) = 0$ on $\Omega$. Then, for $\kappa\in \Omega$, we have
\begin{align}
f(0)&=A^{-1}B, \qquad \nonumber \\
(I+A^{-1}BK&+\kappa A^{-1})f(\kappa)=0 \label{vpinv0}.
\end{align}
In this case, consider the pointwise power series expansion of $f$ around $0$ given by
\begin{equation}\label{expf}
f(\kappa )=A^{-1}B
+\sum_{k=1}^\infty  f_k\kappa^k,
\end{equation}
for some sequence of functions $(f_k)_{k\in\mathbb N^*}$.
Notice that since $f\in D(A+BK)$ and $A^{-1}B\in D(A+BK)$, we obtain that $\sum_{k=1}^\infty  f_k\kappa^k \in D(A)$. By dividing by $\kappa^p$ for any $p\in \mathbb N^*$ and making $\kappa\rightarrow 0$, we deduce that $f_k\in D(A)$ for any $k\in \mathbb N^*$.
At the zeroth order, \eqref{vpinv0} writes
\[
A^{-1}B+A^{-1}BKA^{-1}B=0,
\]
since we assumed that $K(A^{-1}B)=-1$. At the higher orders,  for any $k\in \mathbb N^*$, \eqref{vpinv0} gives
\begin{equation}\label{vpordre}
f_k+A^{-1}BK f_k+\ A^{-1} f_{k-1} =0,
\end{equation}
where $f_0=A^{-1}B$ according to \eqref{expf}. Applying $K$ to \eqref{vpordre} and using $K(A^{-1}B)=-1$ implies that 

\[
K\left(A^{-1}f_k\right)=0, \quad \forall k\geq 0.
\]
By successively taking $A^{-1}$ and $K$ of \eqref{vpordre}, we obtain
\begin{align}\label{Kexp}
K\left(A^{-n}f_k \right)=0, \quad \forall k\geq 0, \quad \forall n\geq 1.
\end{align}
Notably, for $k=0$, we obtain that 
\begin{equation}\label{mar}K\left(A^{-n-1}B\right)=0, \quad \forall n\geq 1.\end{equation}
Remind that $A^{-1}B(x)=x^{1-\alpha}$ has already been decomposed as $A^{-1}B= \sum_{k\in \mathbb N^*}b_k \phi_k$, where the $b_k$ are given in \eqref{Abnd}. Then, using \eqref{Kfd}, we deduce that 
\begin{align}\label{K1}
\sum_{j\in \mathbb N^*} (-1)^j\psi_j(1)\frac{b_j}{\lambda_j^n}=0, \quad \forall n\geq 1.
\end{align}

Consider the entire function
\[
H(z):=\sum_{j\in \NN^*} \frac{ -\psi_j(1)b_j e^{-z/\lambda_j }}{\lambda_j}.
\]
From \eqref{K1}, we obtain that for any $p\in\mathbb N$,  $H^{(p)}(0)=0$ and therefore $H\equiv 0$. By letting $z \rightarrow -\infty$ and using that  
$b_j\not = 0$ for all $j\in \mathbb N^*$ by \eqref{Abnd} and \eqref{behaviorB}, we deduce that 
$\psi_j(1)=0$ for all $j\in \mathbb N^*$.
This means that $\psi(1)=0$. However, this is impossible because of \eqref{pp}, \eqref{wtpsin} and \eqref{nonresonance}.
Hence, we have a a contradiction with the fact that $\lambda(\kappa)=0$ on $\Omega$. Hence, there exists at least one $\kappa\in \Omega$ such that $(I+A^{-1}BK+\kappa A^{-1})$ is invertible. Setting $\rho=\kappa-\lambda$, we deduce that $(A+BK+(\lambda + \rho) I)$ is invertible.

In both cases ($K(A^{-1}B)\not=-1$ and $K(A^{-1}B)=-1$), there exists $\rho\in \mathbb C$ such that  $(A+BK+(\lambda + \rho) I)$ is invertible, this property being still true in a small neighbourhood $\mathcal O$ of $\rho$. As a consequence, since $A+\rho I$ has discrete eigenvalues, it is possible  to modify $\rho$ (by choosing another point of $\mathcal O$) in such a way that $A+\rho I$ is also invertible, which concludes the proof.

\cqfd

\small{
\bibliographystyle{plain}
\bibliography{biblio}}
%
%
%
%
%
%

\end{document}